\documentclass[10pt,draft]{amsart}
\usepackage{amsmath,amssymb,amsthm}
\begin{document}
\newtheorem{lem}{Lemma}[section]
\newtheorem{prop}{Proposition}[section]
\newtheorem{cor}{Corollary}[section]
\numberwithin{equation}{section}
\newtheorem{thm}{Theorem}[section]
\theoremstyle{remark}
\newtheorem{example}{Example}[section]
\newtheorem*{ack}{Acknowledgment}
\theoremstyle{definition}
\newtheorem{definition}{Definition}[section]
\theoremstyle{remark}
\newtheorem*{notation}{Notation}
\theoremstyle{remark}
\newtheorem{remark}{Remark}[section]
\newenvironment{Abstract}
{\begin{center}\textbf{\footnotesize{Abstract}}%
\end{center} \begin{quote}\begin{footnotesize}}
{\end{footnotesize}\end{quote}\bigskip}
\newenvironment{nome}
{\begin{center}\textbf{{}}%
\end{center} \begin{quote}\end{quote}\bigskip}

\newcommand{\triple}[1]{{|\!|\!|#1|\!|\!|}}
\newcommand{\xx}{\langle x\rangle}
\newcommand{\ep}{\varepsilon}
\newcommand{\al}{\alpha}
\newcommand{\be}{\beta}
\newcommand{\de}{\partial}
\newcommand{\la}{\lambda}
\newcommand{\La}{\Lambda}
\newcommand{\ga}{\gamma}
\newcommand{\del}{\delta}
\newcommand{\Del}{\Delta}
\newcommand{\sig}{\sigma}
\newcommand{\ome}{\omega}
\newcommand{\Ome}{\Omega}
\newcommand{\C}{{\mathbb C}}
\newcommand{\N}{{\mathbb N}}
\newcommand{\Z}{{\mathbb Z}}
\newcommand{\R}{{\mathbb R}}
\newcommand{\Rn}{{\mathbb R}^{n}}
\newcommand{\Rnu}{{\mathbb R}^{n+1}_{+}}
\newcommand{\Cn}{{\mathbb C}^{n}}
\newcommand{\spt}{\,\mathrm{supp}\,}
\newcommand{\Lin}{\mathcal{L}}
\newcommand{\SSS}{\mathcal{S}}
\newcommand{\F}{\mathcal{F}}
\newcommand{\xxi}{\langle\xi\rangle}
\newcommand{\eei}{\langle\eta\rangle}
\newcommand{\xei}{\langle\xi-\eta\rangle}
\newcommand{\yy}{\langle y\rangle}
\newcommand{\dint}{\int\!\!\int}
\newcommand{\hatp}{\widehat\psi}
\renewcommand{\Re}{\;\mathrm{Re}\;}
\renewcommand{\Im}{\;\mathrm{Im}\;}

\title[On the Local Smoothing]%
{{On the Local Smoothing for a class of conformally invariant 
Schr\"odinger equations}}
\author{}

\author{Luis Vega}

\address{Luis Vega\\
Universidad del Pais Vasco, Apdo. 64\\
48080 Bilbao, Spain}

\email{mtpvegol@lg.ehu.es}

\author{Nicola Visciglia}

\address{Nicola Visciglia\\
Dipartimento di Matematica Universit\`a di Pisa\\
Largo B. Pontecorvo 5, 56100 Pisa, Italy}

\email{viscigli@dm.unipi.it}

\thanks{\noindent This research was supported by 
RTN Harmonic Analysis and Related Problems, contract HPRN-CT-2001-00273-HARP.
The second author 
was supported also 
by the INDAM project "Mathematical Modeling and Numerical
Analysis of Quantum Systems
with Applications to Nanosciences"}

\maketitle
\date{}
\begin{Abstract}
We present some a - priori bounds from above and 
from below for solutions to a class of conformally invariant 
Schr\"odinger equations. As a by - product we deduce some new uniqueness
results.
\end{Abstract}

\section{Introduction}

The main aim of this paper is to extend a previous result 
proved in \cite{VV} about the local smoothing 
for the free Schr\"odinger equation,
to a more general class of Schr\"odinger 
equations (linear and semilinear)
that are invariant under the conformal transformation.

More precisely we shall consider the following Cauchy problems:
\begin{equation}\label{SL}
{\bf i} \partial_t u - \Delta u + |x|^{-2} 
W\left (\frac{x}{|x|}\right ) u=0,
\end{equation}
$$
u(0)=f, (t, x) \in {\mathbf R}\times {\mathbf R}^n, n\geq 3,
$$
where $W: {\mathbf S}^{n-1}\rightarrow \mathbf R$ is a non - negative,
bounded and measurable function (see also remark \ref{smallBurq}) and
\begin{equation}\label{SN}
{\bf i} \partial_t u - \Delta u \pm u|u|^\frac 4n=0,
 \end{equation}
$$u(0)=f, (t, x) \in {\mathbf R}\times {\mathbf R}^n, n\geq 3.$$

The main property shared by the Cauchy problems 
\eqref{SL} and \eqref{SN} is that both are invariant
under the conformal transformation.
Let us recall that the conformal transformation is the map
$$u(t, x)\rightarrow \tilde u(t, x),$$
defined as follows:
\begin{equation}
\label{lec}\tilde u(t, x)
=\frac{1}{t^\frac n2} e^{\frac{{\bf i}|x|^2}{4t}} 
u \left (\frac 1t, \frac xt \right),
(t, x)\in (0, \infty) \times {\mathbf R}^n.
\end{equation}

This transformation has been extensively used in the literature
in connection with the Schr\"odinger equation 
(for more details see \cite{caze} and the bibliography therein).
In fact an explicit computation shows that if $u(t, x)$ satisfies
\eqref{SL} and $\tilde u(t, x)$ is defined as in \eqref{lec}, then  
\begin{equation}\label{SLP}
{\bf i} \partial_t \tilde u + \Delta \tilde u
- |x|^{-2} W\left (\frac{x}{|x|}        \right ) \tilde u=0,
\end{equation}
\begin{equation*}
\tilde u(1)= e^{{\bf i}\frac{|x|^2}4} u(1), 
(t, x)\in (0, \infty) \times {\mathbf R}^n, n\geq 3.
\end{equation*}

Similarly if $u(t, x)$ satisfies  
\eqref{SN}, then the corresponding $\tilde u(t, x)$
satisfies:
\begin{equation}\label{SNP}
{\bf i} \partial_t \tilde u + \Delta \tilde u \mp 
\tilde u|\tilde u|^\frac 4n=0,
\end{equation}
\begin{equation*}
\tilde u(1)= e^{{\bf i}\frac{|x|^2}4} u(1),
(t, x)\in (0, \infty) \times {\mathbf R}^n, n\geq 3.
\end{equation*}

As it has been mentioned above this article is mainly devoted
to study the local smoothing
for the solutions to the Cauchy problems \eqref{SL}
cand \eqref{SN}, see \cite{CS}, \cite{Sj}, \cite{Vega}. More 
precisely we shall present
some estimates, from above and from below, related
with the phenomena of gain
of $\frac 12$ - derivative for the solution to \eqref{SL} and 
\eqref{SN}.

Let us recall also that in the free case
(i.e. the linear Schr\"odinger equation with
constant coefficient) these estimates have been already proved in 
\cite{VV}.
We were motivated by the results of Agmon and H\"ormander in \cite{AH}.
As it will be clear in the sequel the results that we shall prove
will allow us to deduce some new uniqueness criteria for solutions
to the Cauchy problems \eqref{SL} and \eqref{SN}.

\vspace{0.1cm}

The first result that we shall present 
concerns the Cauchy problem \eqref{SL}.

\begin{thm}\label{mainSL}
Assume that $W:{\mathbf R}^n\rightarrow \mathbf R$ is bounded, measurable
and non - negative and let $u$ be the unique solution to \eqref{SL}
with initial data $f\in \dot H^\frac 12({\mathbf R}^n)$,
then the following a priori estimate is satisfied:
\begin{equation}\label{eqSL}
c \|f\|_{\dot H^\frac 12({\mathbf R}^n)}^2 
\leq \sup_{R>0}\frac 1R \int_{0}^\infty \int_{|x|<R}
|\nabla_x u|^2 \hbox{ } dxdt \leq C \|f\|_{\dot H^\frac 12({\mathbf R}^n)}^2
\end{equation}
where $c, C>0$ are suitable constants.
\end{thm}

\begin{remark}
As it will be clear in the sequel, 
the proof of \eqref{eqSL} will be done by a density argument.
Hence we can assume that the initial data $f$ is 
regular enough in order to guarantee
the existence and the uniquess of solution to \eqref{SL}.
\end{remark}
\begin{remark}
Let us point - out
that the r.h.s. estimate in \eqref{eqSL} has 
been proved already in the paper
\cite{BRV}. Then our main contribution 
is the l.h.s. estimate in \eqref{eqSL}. 
Let us recall also that in \cite{VV} the estimate \eqref{eqSL}
has been proved for the free Schr\"odinger equation, 
i.e. \eqref{SL} with $W\equiv 0$.
\end{remark}

In fact the l.h.s. 
in \eqref{eqSL} will follow from the following

\begin{thm}\label{mainUCSL}
Let $W$, $u$ and $f$ as in theorem \ref{mainSL}.
Then there exists a constant $c>0$ such that
\begin{equation}\label{eqUCSL}
\liminf_{R\rightarrow \infty} \frac 1R 
\int_{0}^\infty \int_{|x|<R}|\nabla_x u|^2 \hbox{ } dxdt
\geq c \|f\|_{\dot H^\frac 12({\mathbf R}^n)}^2 
\hbox{ } \forall f\in \dot H^\frac 12({\mathbf R}^n).
\end{equation}
Therefore if 
$$\liminf_{R\rightarrow \infty} \frac 1R \int_{0}^\infty 
\int_{|x|<R}|\nabla_x u|^2 \hbox{ } dxdt
=0$$
then $u\equiv0$.
\end{thm}

\begin{remark}
Let us recall that the Cauchy problem \eqref{SL}
has been  
investigated in \cite{bpst} from the point of view of
Strichartz estimates, while  
in \cite{BRV} it has been studied 
in connection with
the local smoothing phenomena. 
\end{remark}

\begin{remark}\label{smallBurq}
The non - negativity assumption done on $W$ 
in theorems \ref{mainSL} and \ref{mainUCSL}
could be relaxed by assuming
a smallness condition on its negative part. 
However for simplicity we assume $W\geq 0$ in order
to avoid technical difficulties and to make more transparent
the idea of the proof.
\end{remark}

As a by - product of the techniques involved in the 
proof of the previous theorems, 
we can get another uniqueness result for solutions to \eqref{SL}.
As far as we know the content of next result it 
is not explicitely written elsewhere also in the 
case of the free Schr\"odinger equation.
\begin{thm}\label{uniqueSL}
Let $u$ satisfies \eqref{SL}, with $W\geq 0$ and 
bounded and $f\in L^2({\mathbf R}^n)$.
If you assume that $$\liminf_{t\rightarrow \infty} 
\frac 1t \int_{\mathbf R^n}
|x| |u(t, x)|^2 dx=0,$$
then $u\equiv 0$.
\end{thm}

Next we shall present the corresponding version
of the previous theorems 
for the solutions to
\eqref{SN}.

\begin{thm}\label{mainSN}
Let $\epsilon_0>0$ be a small parameter such that \eqref{SN} has a 
unique global solution $u\in {\mathcal C}({\mathbf R}; H^1({\mathbf R}^n)
\cap L^{2+\frac 4n}({\mathbf R}
\times {\mathbf R}^n)$, for any $f$ such that $\|f\|_{L^2({\mathbf R}^n)}<\epsilon_0$.
There exists $0<\epsilon\leq \epsilon_0$
such that if $u$ is the unique solution to \eqref{SN} 
with initial data $f$ that satisfies 
$\|f\|_{L^2({\mathbf R}^n)}<\epsilon$ and moreover
$$f\in H^1({\mathbf R}^n) \hbox{ and } \int_{{\mathbf R}^n} 
|x|^2|f(x)|^2 dx<\infty,$$
then the following estimate holds: 
\begin{equation}\label{eqSN}
c \|f\|_{\dot H^\frac 12({\mathbf R}^n)}^2 
\leq \sup_{R>R_0}\frac 1R \int_{0}^\infty \int_{|x|<R}
|\nabla_x u|^2 \hbox{ } dxdt \end{equation}
$$\leq C \left( \|f\|_{\dot H^\frac 12({\mathbf R}^n)}^2 
+ \frac 1{R_0^{1-\frac 2n}}
\|f\|_{\dot H^\frac 12({\mathbf R}^n)}^\frac 4n
\right ) \hbox{ } \forall R_0>0$$
where $c, C>0$ are universal constants independent of $f$ and $R_0>0$.
Moreover we have the following chain of inequalities:
\begin{equation}\label{eqUCSN}
c \|f\|_{\dot H^\frac 12({\mathbf R}^n)}^2 
\leq \liminf_{R\rightarrow \infty}\frac 1R \int_{0}^\infty \int_{|x|<R}
|\nabla_x u|^2 \hbox{ } dxdt \end{equation}
$$\leq \limsup_{R\rightarrow \infty}
\frac 1R \int_{0}^\infty \int_{|x|<R}
|\nabla_x u|^2 \hbox{ } dxdt
\leq C \|f\|_{\dot H^\frac 12({\mathbf R}^n)}^2.$$
\end{thm}

\begin{remark}
Let us point - out that the statement of
theorem \ref{mainSN} contains a global existence result
for the Cauchy problem \eqref{SN} under a smallness 
assumption on the initial data $f$
in $L^2({\mathbf R}^n)$.
Let us recall that this fact has been proved in \cite{cw}, \cite{cw2} 
and \cite{t}. 
Hence our main contribution
concerning the Cauchy problem \eqref{SN}
are the estimates \eqref{eqSN} and \eqref{eqUCSN}. 
\end{remark}

We shall prove also the following nonlinear version of theorem 
\ref{uniqueSL}.

\begin{thm}\label{uniqueSN}
Let $u$ be the unique global solution to \eqref{SN}
with $\|f\|_{L^2({\mathbf R}^n)}<\epsilon_0$, where $\epsilon_0>0$
is small enough.
Assume moreover that 
$$\liminf_{t\rightarrow \infty} \frac 1t \int_{\mathbf R^n}
|x| |u(t, x)|^2 dx=0,$$
then $u\equiv 0$.
\end{thm}

Along the proof of theorem \ref{mainSN} (in particular 
in the proof of the l.h.s. estimate
in \eqref{eqUCSN})
we shall need some intermediate results that in our 
opinion have their own interest.
One of them will be stated in next
theorem.
\begin{thm}\label{radiation}
Let $\psi\in C^1({\mathbf R}^n)$ be a 
radially symmetric function such that the following limit exists:
\begin{equation}
\label{infbil}\lim_{|x|\rightarrow \infty} \partial_{|x|}\psi(x)
=\psi'(\infty)\in [0, \infty),
\end{equation}
and moreover $$\partial_{|x|}\psi(x)\geq 0
\hbox{   } \forall x\in {\mathbf R}^n.$$
Assume that $u$ and $f$  
are as in theorem \ref{mainSN}, 
then the following estimate holds:
\begin{equation}\label{somcon7NLS}
\limsup_{t\rightarrow \infty} \left (
- {\mathcal Im}\int_{{\mathbf R}^n} \bar u(t) \nabla u(t) 
\cdot \nabla \psi(x)  \hbox{ } 
dx
\right )\geq \frac 12 \psi'(\infty) \int_{{\mathbf R}^n} 
|x||g(x)|^2 \hbox{ } dx,
\end{equation}
where $g$ is a suitable function 
that depends on $f$ but does not depend on $\psi$.
Moreover the following inequality is satisfied:
\begin{equation}\label{dotNLS}
\|f\|_{\dot H^\frac 12 ({\mathbf R}^n)}^2 \leq C 
\int_{{\mathbf R}^n} |x||g(x)|^2 \hbox{ } dx,
\end{equation}
where $C>0$ is suitable constant that does not depend on $f$.
\end{thm}

\begin{remark}
Looking at the proof of theorem \ref{radiation} it will be clear that
the smallness assumption done on the initial condition $f$
in $L^2({\mathbf R}^n)$ will be relevant in order to prove
\eqref{dotNLS}. However in the defocusing case 
(i.e. \eqref{SN} with the sign plus),
the existence of the function $g$,
the validity of \eqref{somcon7NLS} and theorem \ref{uniqueSN} can be proved without any
smallness assumption on $f$ and just assuming $f\in H^1({\mathbf R}^n)$.
\end{remark}

\begin{remark}
In the proof of theorems \ref{mainSL} and 
\ref{mainUCSL} we shall need
a result similar to theorem \ref{radiation} concerning 
the solutions to \eqref{SL}
(see proposition  \ref{radiationSL}). 
However the proof of theorem \ref{radiation} (and in particular
the proof of \eqref{dotNLS}) 
is much more involved due to the 
nonlinear nature 
of the operator $f\rightarrow g$ given in the statement 
of theorem \ref{radiation}. 
\end{remark}

Next we shall fix some notations that will be used 
in the sequel.

\vspace{0.1cm}

{\bf Notations.}
For any $s\in \mathbf R$ we shall 
denote by $\dot H^s_x$ and $H^s_x$ the
homogeneous and non - homogeneous Sobolev 
spaces in ${\mathbf R}^n$ of order $s$.

For any $R>0$ we shall denote by $B_R$
the unit ball in ${\mathbf R}^n$
centered in the origin.

For any $1\leq p, q\leq \infty$ 
$$L^p_x \hbox{ and } L^p_t L^q_x$$
denote the Banach spaces
$$L^p({\mathbf R}^n) \hbox{ and } L^p({\mathbf R}; L^q({\mathbf R}^n)).$$
We shall also write
$$L^p_t L^p_x= L^p_{t, x}.$$

Let $X$ be a general Banach spaces,
then ${\mathcal C}_t (X)$
is the space of continuous functions defined in $\mathbf R$
and valued in $X$. 

Given any non - negative and measurable function 
$w:{\mathbf R}^n\rightarrow {\mathbf R}^+$
we shall denote by $L^2_w$ the Hilbert space
whose norm is defined as follows:
$$\|f\|_{L^2_w}^2=\int_{{\mathbf R}^n} |f(x)|^2 w(x) dx.$$

Given a space - time dependent function $w(t, x)$ we shall denote by
$w(t_0)$ the trace of $w$ at fixed 
time $t\equiv t_0$, in case that it is 
well - defined.

We shall denote by $\int ...\hbox {  }dx ,
\int ... \hbox{ } dt$
and $\int \int ...\hbox{ } dx dt$ the 
integral of suitable functions
with respect to the full space, time  
and space - time variables respectively.

When it is not better specified we shall 
denote by $\nabla v$ the gradient
of any time - dependent function $v(t, x)$ 
with respect to the space variables.
Moreover $\nabla_\tau$ and $\partial_{|x|}$
shall denote respectively the tangential 
gradient and the radial derivative.

If $\psi\in C^2({\mathbf R}^n)$, then 
$D^2 \psi$ will represent the hessian matrix of $\psi$.
\section{The conformal conservation law}

In order to simplify the proof of our results 
it will be useful in some cases to work directly
with the following general class of Cauchy problems:
\begin{equation}\label{SM}
{\bf i} \partial_t u - \Delta u + |x|^{-2}W\left 
(\frac x{|x|} \right) u \pm \lambda u|u|^\frac 4n=0, 
\end{equation}
\begin{equation*}
u(0)=f, \lambda \geq 0, (t, x) \in {\mathbf R}\times {\mathbf R}^n
\end{equation*}
and
\begin{equation}\label{SMP}
{\bf i} \partial_t u + \Delta u - |x|^{-2}W
\left( \frac x{|x|}\right ) u \mp \lambda u|u|^\frac 4n=0, 
\end{equation}
\begin{equation*}
u(0)=f, \lambda \geq 0, (t, x) \in {\mathbf R}\times {\mathbf R}^n.
\end{equation*}

The following result can be found in \cite{caze}.

\begin{prop}\label{sommerfield}
Let $u\in {\mathcal C}_t (H^1_x)$ satisfies
\eqref{SM} with $W\geq 0$ and
$f\in H^1_x\cap L^2_{|x|^2}$, then:
\begin{equation}
\label{energy}\int \left (|\nabla u(t)|^2 + |x|^{-2}
W\left(\frac x{|x|} \right) |u(t)|^2 
\pm \frac {\lambda n}{n+2} |u(t)|^{2+\frac 4n} \right )dx
\end{equation}
\begin{equation*}
=\int \left (|\nabla f(x)|^2 + |x|^{-2}
W\left (\frac x{|x|} \right) |f(x)|^2 
\pm \frac {\lambda n}{n+2} |f(x)|^{2+\frac 4n} \right )dx \hbox{  } 
\forall t\in \mathbf R;
\end{equation*}
\begin{equation}\label{charge}
\int |u(t)|^2 dx = \int |f(x)|^2 dx \hbox{  } \forall t\in \mathbf R;
\end{equation}
\begin{equation}\label{conservation}
\|x u(t) - 2{\bf i}t \nabla u(t)\|_{L^2_x}^2 + 
 4 t^2 \int |x|^{-2}W\left (\frac x{|x|}\right ) |u(t)|^2 dx 
\end{equation}
\begin{equation*}
\pm  \frac{2n \lambda t^2 }{n+2}  \int |u(t)|^{2+\frac 4n} dx 
= \|f\|_{L^2_{|x|^2}}^2 \hbox{  } \forall t\in \mathbf R.
\end{equation*}
\end{prop}

Notice that if we choose in a suitable way the parameter $\lambda$ 
and the potential $W$ in \eqref{SM}, then we 
can deduce from the previous proposition
the following
corollary.
\begin{cor}
Assume that $u \in {\mathcal C}_t(H^1_x)$ 
solves \eqref{SL} with $W\geq 0$
and $f \in H^1_x\cap L^2_{|x|^2}$,
then:
\begin{equation}\label{conservation1}
\|x u(t) - 2{\bf i}t \nabla u(t)\|_{L^2_x}^2 + 
 4 t^2 \int |x|^{-2}W\left (\frac x{|x|} \right) |u(t)|^2 dx  
= \|f\|_{L^2_{|x|^2}}^2 \hbox{ } \forall t \in \mathbf R.
\end{equation}
If $u\in  {\mathcal C}_t(H^1_x)$ solves \eqref{SN} with 
$f \in H^1_x\cap L^2_{|x|^2}$, then:
\begin{equation}\label{conservation2}
\|x u(t) - 2{\bf i}t \nabla u(t)\|_{L^2_x}^2 
\pm  \frac{2n t^2 }{n+2}  \int |u(t)|^{2+\frac 4n} dx 
= \|f\|_{L^2_{|x|^2}}^2 \hbox{ } \forall t \in \mathbf R.
\end{equation}
\end{cor}

The following proposition is similar to proposition 
\ref{sommerfield} except that it concerns the solutions
to \eqref{SMP}.

\begin{prop}\label{sommerfieldP}
Assume that $u\in {\mathcal C}_t(H^1_x)$ satisfies \eqref{SMP}
where $f\in L^2_x \cap L^2_{|x|^2}$, 
then the following 
identities are satisfied:
\begin{equation}\label{charge2}
\int |u(t)|^2 dx = \int |f(x)|^2 dx \hbox{  } \forall t\in \mathbf R;
\end{equation}

\begin{equation}\label{conservation3}
\|x u(t) + 2{\bf i}t \nabla u(t)\|_{L^2_x}^2 + 
 4 t^2 \int |x|^{-2}W\left (\frac x{|x|} \right ) |u(t)|^2 dx 
\end{equation}
\begin{equation*}
\pm  \frac{2n \lambda t^2 }{n+2}  \int |u(t)|^{2+\frac 4n} dx 
= \|f\|_{L^2_{|x|^2}}^2
\hbox{  } \forall t\in \mathbf R. 
\end{equation*}
\end{prop}

Notice that as a by - product of proposition \ref{sommerfieldP} 
(where we choose $W\equiv 0$ and $\lambda=0$) we get the following

\begin{cor}\label{pseudoNLS}
Assume that
$${\bf i}\partial_t u + \Delta u=0, $$$$u(0)=f\in L^2_{|x|}$$
then
$$\|e^{-{\bf i} \frac{|x|^2}{4t}}u(t)
\|_{\dot H^\frac 12_x}^2 \leq \frac 1{2t} \|f \|_{L^2_{|x|}}^2\hbox{  }
\forall t\in (0, \infty).$$
\end{cor}

\noindent {\bf Proof.}
By using  
\eqref{conservation3}, where we choose $W\equiv 0$ and $\lambda=0$,
we get:
$$4t^2\|e^{-{\bf i}\frac{|x|^2}{4t}}u(t)\|_{\dot H^1_x}^2
=
4 t^2 \|\nabla (e^{-{\bf i}\frac{|x|^2}{4t}}u(t))\|_{L^2_x}^2
= \|f \|_{L^2_{|x|^2}}^2.$$

On the other hand the conservation of the charge (see \eqref{charge2})
implies
$$\|e^{-{\bf i}\frac{|x|^2}{4t}}u(t)\|_{L^2_x}^2=\|u(t)\|_{L^2_x}
=\|f \|_{L^2_x}^2.$$

The result follows
by interpolation.
 
\hfill$\Box$

\section{On the asymptotic behaviour of solutions to \eqref{SL}\\
and proof of theorem \ref{uniqueSL}}

The main result of this section is the following

\begin{prop}\label{radiationSL} 
Let $\psi\in C^1({\mathbf R}^n)$ be a 
radially symmetric function such that the following limit exists:
\begin{equation}\label{limitder}\lim_{|x|\rightarrow \infty} 
\partial_{|x|}\psi=\psi'(\infty)\in [0, \infty),
\end{equation}
and moreover 
\begin{equation}\label{nonneg}\partial_{|x|}\psi\geq 0
\hbox{  } \forall x\in {\mathbf R}^n.
\end{equation}
Let $u\in {\mathcal C}_t(H^1_x)$ be the 
unique global solution to \eqref{SL},
where $W$ is bounded and non - negative 
and $f\in H^1_x \cap L^2_{|x|^2}$, 
then 
\begin{equation}\label{somcon7}
\liminf_{t\rightarrow \infty} \left (
- {\mathcal Im}\int \bar u(t) \nabla u(t) 
\cdot \nabla \psi  \hbox{ } dx
\right )\geq \frac 12 \psi'(\infty) \int |x||g(x)|^2 dx,
\end{equation}
where $g$ is a 
suitable function that depends on $f$ but does not depend on $\psi$.
Moreover the following estimate holds:
\begin{equation}\label{dot}
\|f \|_{\dot H^\frac 12_x}^2 \leq \frac 12 \int |x||g(x)|^2 dx.
\end{equation}
\end{prop}

\noindent {\bf Proof.}
Let us notice that \eqref{conservation1} implies the following identity:
\begin{equation}\label{pseudoestimateSL}
\left \|\frac xt u(t) - 2{\bf i}\nabla u(t)\right \|_{L^2_x}^2 + 
 4 \int |x|^{-2}W\left (\frac x{|x|}\right) |u(t)|^2 dx 
= \frac 1{t^2}\|f \|_{L_{|x|^2}^2}^2,
\end{equation}
that, due to the non - negativity assumption done
on $W$, implies:
\begin{equation}\label{som}
\lim_{t\rightarrow \infty}\left \|\frac xt u(t) - 2{\bf i} \nabla u(t)
\right \|_{L^2_x}^2
=0.
\end{equation}

We split now the proof in two parts.
\\
\\
{\em Construction of $g$ and proof of \eqref{dot}}
\\
\\

Let us recall that if
$u$ satisfies \eqref{SL}, then its conformal transformation
$$\tilde u (t, x)= \frac{1}{t^\frac n2} e^{\frac{{\bf i}|x|^2}{4t}} 
u \left (\frac 1t, \frac xt \right),
$$
satisfies the Cauchy problem \eqref{SLP}.
In particular 
\begin{equation}\label{trivPseu}\tilde u(1)
= e^{\frac{{\bf i}|x|^2}{4}} u \left (1\right)
\end{equation}
and hence 
$$\|\tilde u(1)\|_{L^2_x}=\|u(1)\|_{L^2_x}= \|f \|_{L^2_x},$$
where we have used the conservation 
of the charge for the unique solution to
\eqref{SL} (see \eqref{charge}).
As a consequence $\tilde u(t, x)$ satisfies 
the following Cauchy problem 
\begin{equation}\label{idido}
{\bf i}\partial_t \tilde u 
+ \Delta \tilde u - |x|^{-2}W\left (\frac x{|x|} \right)\tilde u=0, 
\end{equation}
\begin{equation*}
\tilde u(1)\in L^2_x, (t, x)\in (0,\infty)\times {\mathbf R}^n.
\end{equation*}

Due to the global well - posedness (in the $L^2_x$ sense) 
of the previous Cauchy problem, 
we deduce that $\tilde u$ can be extended
as a solution to the same Cauchy problem in
the functional space ${\mathcal C}_t( L^2_x)$.
In particular it is
well defined a function
$g \in L^2_x$ such that
the following limit exists in $L^2_x$:
\begin{equation}\label{lim}
\lim_{t\rightarrow 0} \tilde u (t, x)= g \in L^2_x. 
\end{equation}

Due to
\eqref{lim} we can deduce that
$$\lim_{t\rightarrow 0}
\left \| u\left ( \frac 1t, x\right)- 
t^\frac n2 e^{-{\bf i} t \frac{|x|^2}{4}} 
g \left(tx\right)\right \|_{L^2_x}$$
$$=\lim_{t\rightarrow 0}
\left \|\frac 1{t^\frac n2} e^{{\bf i} \frac{|x|^2}{4t}} 
u\left (\frac 1t, \frac xt \right)
- g \left (x\right )\right \|_{L^2_x}=0,$$
and in particular
\begin{equation}\label{L2SL}
\lim_{t\rightarrow \infty} \left \|u(t) - 
\frac 1{t^\frac n2}e^{-{\bf i} \frac{|x|^2}{4t}} g\left
  (\frac{x}{t}\right )\right \|_{L^2_x}=0.
\end{equation}

On the other hand \eqref{idido} and \eqref{lim} imply that 
$\tilde u$ satisfies
$${\bf i}\partial_t \tilde u + \Delta \tilde u 
- |x|^{-2}W\left ( \frac x{|x|}\right)
\tilde u =0$$$$ \tilde u (0)=g$$
that in turn, due to \eqref{conservation3}
(where we choose $\lambda=0$ and $t=1$),
implies
$$\|x\tilde u (1) + 2{\bf i}\nabla \tilde u (1)\|_{L^2_x}^2 +
4\int |x|^{-2}W\left( \frac x{|x|}\right)
|\tilde u (1)|^2 dx= \|g\|^2_{L^2_{|x|^2}}.$$

Notice that this identity is equivalent to 
$$ 4\|\nabla (e^{-{\bf i}\frac{|x|^2}4} \tilde u(1))\|_{L^2_x}^2 
+
4\int |x|^{-2}W\left( \frac x{|x|}\right)|
e^{-{\bf i}\frac{|x|^2}4} \tilde u (1)|^2 dx = \|g\|_{L^2_{|x|^2}}^2,$$
that due to \eqref{trivPseu}
gives:
\begin{equation}\label{spect}4 \|\nabla u(1)\|_{L^2_x}^2 + 
4\int |x|^{-2}W\left( \frac x{|x|}\right)|u(1)|^2 dx
=\|g \|_{L^2_{|x|^2}}^2.\end{equation}

By combining this identity with \eqref{energy} 
(where we choose $\lambda=0$) and
with the non - negativity assumption done on $W$ we get: 
\begin{equation}\label{subiniti0}
4 \|\nabla f\|_{L^2_x}^2 \leq 4 \|\nabla f\|_{L^2_x}^2 + 
4\int |x|^{-2}W\left( \frac x{|x|}\right)|f(x)|^2 dx
\end{equation}
$$
= 4 \|\nabla u(1)\|_{L^2_x}^2 + 
4\int |x|^{-2}W\left( \frac x{|x|}\right)|u(1)|^2 dx
=\|g \|_{L^2_{|x|^2}}^2.
$$

Notice also that the following estimate follows easily
from the conservation of the charge 
(see \eqref{charge} and \eqref{charge2}):
\begin{equation}\label{chae}\|f \|_{L^2_x}^2 = \|u(1)\|_{L^2_x}^2= 
\|\tilde u(1)\|_{L^2_x}^2=\|g \|_{L^2_x}^2.
\end{equation}

Then \eqref{dot} 
follows by making interpolation 
between \eqref{subiniti0} and \eqref{chae}.
\\  
\\
{\em Proof of \eqref{somcon7}}
\\
\\

Due to \eqref{som} we can deduce that
\begin{equation}\label{som2}\lim_{t\rightarrow \infty}
\left [\int \bar u(t) \nabla u(t) \cdot \nabla \psi \hbox{ } dx 
+\frac {\bf i} {2t} \int |x| | u(t) |^2  
\partial_{|x|}\psi  \hbox{ } dx\right]=0
\end{equation}
and then
\begin{equation}\label{truncSL}\liminf_{t\rightarrow \infty}
\left( -{\mathcal Im}
\int \bar u(t) \nabla u(t) \cdot \nabla \psi dx \right)
\end{equation}
$$
= \frac 12 \liminf_{t\rightarrow \infty} 
\int |x| | u(t) |^2  \partial_{|x|}\psi \hbox{ } \frac{dx}t.
$$

Next we fix a real number $R>0$ and we notice that due to
the non - negativity assumption done on $\partial_{|x|}\psi$ 
(see \eqref{nonneg}) we get:
\begin{equation}\label{limintSL} \int 
|x|| u(t) |^2 \partial_{|x|}\psi \hbox{ } \frac {dx}t
\geq \int_{|x|\leq Rt} |x| | u(t) |^2 
\partial_{|x|}\psi  \hbox{ } \frac{dx}t
\end{equation}
$$=  \int_{|x|\leq Rt} |x| \left (|u(t)|^2 -\frac 1{t^n} 
\left | g\left (\frac xt \right )\right |^2\right )
\partial_{|x|}\psi \hbox{ }\frac{dx}t
$$$$+ \int_{|x|\leq Rt} |x| \left ((\partial_{|x|} \psi 
-\psi'(\infty)\right)  \left | g
\left ( \frac x {t}\right)\right |^2 \frac{dx}{t^{n+1}}
$$
$$
+ \psi'(\infty) \int_{|x|\leq Rt} 
|x| \left |g\left (\frac xt \right )
\right |^2 \frac{dx}{t^{n+1}} \hbox{ }
\forall R>0,
$$
where $g$ is the function constructed in the previous step.

Notice that the following estimate is trivial:
\begin{equation}\label{lim4SL}
\int_{|x|\leq Rt}
|x| \left (|u(t)|^2 -\frac 1{t^n} 
\left | g\left (\frac xt \right )\right |^2\right )
\partial_{|x|}\psi \hbox{ } \frac{dx}t
\end{equation}
\begin{equation*}
\leq R \|\partial_{|x|}\psi\|_{L^\infty_x}\int_{|x|\leq Rt}
\left ||u(t)|^2 -\frac 1{t^n} 
\left | g\left (\frac xt \right )\right |^2\right |dx
\rightarrow 0 \hbox{ as } t\rightarrow \infty,
\end{equation*}
where at the last step we have used \eqref{L2SL}.

Moreover the change of variable formula implies:
$$\left |\int_{|x|\leq Rt} |x| \left ((\partial_{|x|} \psi 
-\psi'(\infty)\right ) \left | g\left ( \frac x 
{t}\right)\right |^2 \frac{dx}{t^{n+1}}\right |
$$
$$\leq R  \int_{|x|\leq R} \left | \partial_{|x|} \psi ( t x )
-\psi'(\infty)\right| |g(x)|^2 dx,$$ 
that in conjunction  with 
the dominated convergence theorem and with assumption
\eqref{limitder} implies:
\begin{equation}
\label{lim7SL}\lim_{t\rightarrow \infty}\int_{|x|\leq Rt} 
|x| \left ((\partial_{|x|} \psi 
-\psi'(\infty) \right) 
\left | g\left ( \frac x {t}\right)\right |^2 \frac{dx}{t^{n+1}}
=0 \hbox{ } \forall R>0.\end{equation}

Due again to the change of variable formula we get
$$
\psi'(\infty) \int_{|x| \leq Rt} |x| 
\left |g\left (\frac xt \right )\right |^2 \frac{dx}{t^{n+1}}
= \psi'(\infty)\int_{|x|\leq R } |x| |g(x)|^2 dx,
$$
and in particular
\begin{equation}\label{lim14SL}
\lim_{t\rightarrow \infty}
\psi'(\infty) \int_{|x|\leq Rt} 
|x| \left |g\left (\frac xt \right )\right |^2 \frac{dx}{t^{n+1}}
= \psi'(\infty)\int_{|x|\leq R} |x| |g(x)|^2 dx.
\end{equation}

By combining 
\eqref{lim4SL},\eqref{lim7SL}, \eqref{lim14SL} and
\eqref{limintSL} we can deduce that
\begin{equation}\label{limfinSL}
\liminf_{t\rightarrow \infty}\int 
|x| |u(t)|^2 \partial_{|x|}\psi(x)\frac{dx}t
\end{equation}
$$
\geq \psi'(\infty) \int_{|x|\leq R} |x| |g(x)|^2 dx \hbox{  }
\forall R>0.
$$

Since $R>0$ is arbitrary, we can combine 
\eqref{truncSL} with \eqref{limfinSL}
in order to deduce 
\eqref{somcon7}.

\hfill$\Box$

\vspace{0.2cm}

\noindent {\bf Proof of theorem \ref{uniqueSL}}
Let $g$ be the function constructed 
in proposition \ref{radiationSL}.
Looking at the proof of \eqref{limfinSL} it is easy
to deduce with a similar argument the following
estimate:
\begin{equation}\label{uniquenessseb}
\liminf_{t\rightarrow \infty}\int 
|x| |u(t)|^2 \frac{dx}t
\geq \int_{|x| \leq R} |x| |g(x)|^2 dx \hbox{  }
\forall R>0.
\end{equation}

In fact we have:
\begin{equation}\label{limintSL2} \int 
|x|| u(t) |^2 \frac {dx}t
\geq \int_{|x|\leq Rt} |x| | u(t) |^2 \hbox{ } \frac{dx}t
\end{equation}
$$=  \int_{|x|\leq Rt} |x| \left [|u(t)|^2 -\frac 1{t^n} 
\left | g\left (\frac xt \right )\right |^2\right ]
\frac{dx}t
+ \int_{|x|\leq Rt} 
|x| \left |g\left (\frac xt \right )\right |^2 \frac{dx}{t^{n+1}} \hbox{ }
\forall R>0.
$$

\vspace{0.1cm}

Notice that the following estimate is trivial:
\begin{equation}\label{lim4SL2}
\int_{|x|\leq Rt}
|x| \left (|u(t)|^2 -\frac 1{t^n} 
\left | g\left (\frac xt \right )\right |^2\right )\frac{dx}t
\end{equation}
\begin{equation*}
\leq R \int_{|x|\leq Rt}
\left ||u(t)|^2 -\frac 1{t^n} 
\left | g\left (\frac xt \right )\right |^2\right |dx
\rightarrow 0 \hbox{ as } t\rightarrow \infty,
\end{equation*}
where at the last step we have used \eqref{L2SL}.

On the other hand we have:
\begin{equation}\label{lasseb}
\int_{|x| \leq Rt} |x| \left |g\left (\frac xt 
\right )\right |^2 \frac{dx}{t^{n+1}}
= \int_{|x|\leq R } |x| |g(x)|^2 dx.
\end{equation}
and then \eqref{uniquenessseb} follows 
by combining \eqref{limintSL2}, \eqref{lim4SL2},
\eqref{lasseb}.

In particular if 
$$\liminf_{t\rightarrow \infty}\int 
|x| |u(t)|^2 \frac{dx}t=0,$$ then  
\eqref{uniquenessseb} implies  $g\equiv 0$, that in turn due to 
\eqref{L2SL} gives $\lim_{t\rightarrow \infty} \|u(t)\|_{L^2_x}=0$.
By combining this fact with the 
conservation of the charge \eqref{charge},
we get $f\equiv0$, and hence $u\equiv 0$.

\hfill$\Box$

\section{Proof of theorems \ref{mainSL}
and \ref{mainUCSL}}\label{BAR}

In the first part of this section
we recall the approach used in
\cite{BRV} in order to deduce the local smoothing 
estimate (i.e. the r.h.s. in \eqref{eqSL})
for the solutions to \eqref{SL}.
The main idea is to multiply
\eqref{SL} by the quantity 
\begin{equation}\label{multiplier}
\nabla \bar u  \cdot \nabla \psi+ \frac 12 \bar u \hbox{ } \Delta \psi ,
\end{equation}
and to integrate on the strip $(0, T)\times \mathbf R^n$.
For the moment $\psi: {\mathbf R}^n
\rightarrow \mathbf R$ is a general function 
to which we require only minimal regularity assumptions 
in order to justify the integration by parts.

The approach described above allows you to deduce the
following family of identities:
\begin{equation}\label{BarcRuiVeg}
\int_{0}^T\int  \left[\nabla \bar u D^2 \psi \nabla u
-\frac 14  |u|^2\Delta^2 \psi + 
|u|^2 |x|^{-3} W\left( \frac x{|x|}\right) \partial_{|x|}\psi 
\right ]dx dt\end{equation}
\begin{equation*} = -\frac 12
{\mathcal Im}
\int \bar u(T) \nabla u(T) \cdot \nabla \psi \hbox{ } dx + \frac 12
{\mathcal Im}
\int \bar f \hbox{ } \nabla f \cdot \nabla \psi \hbox{  } dx,
\end{equation*}
(for more details on this computation see \cite{BRV} and \cite{VV}).

The following propositions will be relevant in the sequel.

\begin{prop}\label{brvrad}
Let $\psi\in C^1({\mathbf R}^n)$
be a radially symmetric function such that:
\begin{equation*}
|\partial_{|x|}\psi|, |x| |\partial_{|x|}^2\psi|
\leq C <\infty \hbox{  } \forall x \in {\mathbf R}^n,
\end{equation*}
where $C>0$ is a suitable constant.
Then
there exists $C'>0$, that depends only on $C$ and such that:
\begin{equation}
\label{dele23}\left |\int g(x) \nabla \bar g(x) 
\cdot \nabla \psi \hbox{ } dx
\right |\leq C' \|g\|_{\dot H^\frac 12_x}^2 \hbox{  } \forall g \in H^1_x.
\end{equation}
In particular if 
$u\in {\mathcal C}_t(H^1_x)$ is the unique solution to the
Cauchy problem \eqref{SL} where $W\geq 0$ and bounded and $f \in H^1_x$,
then:
\begin{equation}
\label{dele}\left |\int u(t) \nabla \bar u(t) 
\cdot \nabla \psi \hbox{ }  dx
\right |\leq C' \|f\|_{\dot H^\frac 12_x}^2 
\hbox{  } \forall t \in \mathbf R.
\end{equation}
\end{prop}

\noindent {\bf Proof.} The proof of \eqref{dele23}
can be found in \cite{BRV}.
In order to prove \eqref{dele}
let us notice that if we choose $g(x)=u(t)$
in \eqref{dele23}
then the inequality becomes
\begin{equation}
\label{delepar}\left |\int u(t) 
\nabla \bar u(t) \cdot \nabla \psi \hbox{ } dx
\right |\leq C'\|u(t)\|_{\dot H^\frac 12_x}^2.
\end{equation}

On the other hand due to 
\eqref{energy} and by recalling the non - negativity 
assumption done on $W$ we get:
$$\|u(t)\|_{\dot H^1_x}^2\leq 
\|f\|_{\dot H^1_x}^2 + 
\|W\|_{L^\infty_x} \int |x|^{-2}|f(x)|^2 dx\leq
C \|f\|_{\dot H^1_x}^2,$$
where we have used the classical Hardy inequality at the last step.

By making interpolation between this inequality
and $$\|u(t)\|_{L^2_x}^2= \|f\|_{L^2_x}^2$$
that follows from \eqref{charge}, we deduce 
$$\|u(t)\|_{\dot H^\frac 12_x}\leq C 
\|f\|_{\dot H^\frac 12_x} \hbox{ } \forall t\in \mathbf R.$$ 

By combining this last inequality with
\eqref{delepar} we get \eqref{dele}.

\hfill$\Box$

\begin{cor}\label{reM}
Assume that $u\in {\mathcal C}_t(H^1_x)$ solves 
\eqref{SL} with $W\geq0$ and bounded
and $f\in H^1_x$, 
then
\begin{equation}\label{ulTI}
\int \int_{|x|<1} \frac{|u|^2}{|x|^2}
W\left(\frac x{|x|} \right) dx dt<\infty.
\end{equation}
\end{cor}

\noindent{\bf Proof.}
Let us notice that if we choose in \eqref{BarcRuiVeg}
the function $\psi$ to be equal to the function $\phi$
given in proposition \ref{brvpsi} 
(see the Appendix) then we get the following inequality:
\begin{equation}\label{QUaS}
\int \int_{|x|<1}  
|u|^2 |x|^{-3} W\left( \frac x{|x|}\right) 
\partial_{|x|}\phi \hbox{ } dx dt 
\leq C \|f\|_{\dot H^\frac 12_x}^2,
\end{equation}
where we have used also \eqref{dele}.

On the other hand 
by the Taylor formula and by using the properties of $\phi$ we get
$$\partial_{|x|}\phi(|x|)=\partial_{|x|}^2\phi(0)|x| + 0(|x|^2).$$

By combining this identity with
\eqref{QUaS} we get
\begin{equation}\label{qUaS}
\partial_{|x|}^2 \phi(0) \int \int_{|x|<1}   
|u|^2 |x|^{-2} W\left( \frac x{|x|}\right) 
dx dt \end{equation}
$$
\leq C \|f\|_{\dot H^\frac 12_x}^2 + 
C \int \int_{|x|<1}   
|u|^2 |x|^{-1} W\left( \frac x{|x|}\right) 
dx dt. 
$$

On the other hand the H\"older inequality implies:
\begin{equation}\label{QuAs}
\int \int_{|x|<1}  |u|^2 |x|^{-1} W\left( \frac x{|x|}\right) 
dx dt  
\end{equation}
$$\leq \|W\|_{L^\infty_x} \|u\|^2_{L^2_t L^{\frac{2n}{n-2}}_x}
\||x|^{-1}\|_{L^{\frac n2}(|x|<1)}< \infty$$
where we have used in the last step the fact that Strichartz estimates
are satisfied by the solutions to \eqref{SL} (see \cite{bpst}).

Since $\partial^2_{|x|}\phi(0)> 0$ by construction, we
can combine \eqref{qUaS} with \eqref{QuAs} 
in order to get the desired result.

\hfill$\Box$

\begin{remark}
Notice that following \cite{BRV} it is possible to show that
$$\int \int_{|x|<1}\frac{|u|^2}{|x|^{3-\epsilon}}\hbox{ } dxdt<\infty $$ 
for any $\epsilon>0$ and for any $u$ that satisfies \eqref{SL}.
However in order to make this paper self - contained we have 
presented a simplified argument to prove corollary \ref{reM}
that is enough for our purpose.
\end{remark}

Next proposition follows by combining proposition
\ref{radiationSL}
with \eqref{BarcRuiVeg}.

\begin{prop}\label{radiation2}
Assume that $\psi$ satisfies the same assumptions 
as in proposition \ref{radiationSL}
and $u$ satisfies \eqref{SL} with $W\geq0$ 
and bounded and $f\in H^1_x\cap L^2_{|x|^2}$.
Then the following estimate holds:
\begin{equation}\label{BarcRuiVegrad}
\int_0^T \int \left[\nabla \bar u D^2 \psi \nabla u
-\frac 14  |u |^2\Delta^2 \psi + 
|u|^2 |x|^{-3} W\left( \frac x{|x|}\right) \partial_{|x|}\psi 
\right ]dx dt
\end{equation}
$$\geq  \frac 14 \psi'(\infty) \|f\|_{\dot H^\frac 12_x}^2 
+ \frac 12 {\mathcal Im}\int \bar f \hbox{ } 
\nabla f\cdot \nabla \psi \hbox{ } dx.$$
\end{prop}

We shall need also the following
\begin{lem}\label{noma}
Assume that $u\in {\mathcal C}_t(H^1_x)$ 
satisfies \eqref{SL} with $f\in H^1_x$ 
and $W\geq 0$ and bounded,
then:
\begin{equation}\label{uno}
\lim_{R\rightarrow \infty}  \int \int 
|u|^2 |\Delta^2 \phi_R| 
\hbox{ } dxdt =0,\end{equation}
\begin{equation}\label{due}
\lim_{R\rightarrow \infty}  \int \int  
\frac{|u|^2} {|x|^{3}} W\left( \frac x{|x|}\right) 
|\partial_{|x|} \phi_R |\hbox{ } dxdt=0
\end{equation}
where $\phi\in C^4({\mathbf R}^n)$ is a radially symmetric function
such that 
$$\partial_{|x|} \phi(0)=0,$$
$$
|\partial_{|x|}\phi| \leq C, |\Delta^2 \phi |\leq \frac{C}
{(1+|x|)^3} \hbox{  } \forall x\in {\mathbf R}^n$$
and $\phi_R =R \phi\left (\frac{x}R\right )$.
\end{lem}

\noindent{\bf Proof.}
\\
{\em Proof of \eqref{uno}}
\\
\\

Notice that the H\"older inequality implies
$$\int \int |u|^2 |\Delta^2 \phi_R| \hbox{ } dx dt
\leq \left [ \int \left ( \int |u(t)|^
\frac{2n}{n-2} \hbox{ } dx \right )^\frac{n-2}{n} dt
\right ]
\left ( \int |\Delta^2 
\phi_R|^\frac n2  \hbox{ } dx \right )^\frac 2n$$
$$
\leq C  
\|u\|_{L^2_t L^\frac{2n}{n-2}_x}^2  
\left (\int  \frac 1{(R + |x|)^\frac {3n}2}\hbox{  } dx \right)^\frac 2n
\rightarrow 0 \hbox{ as
} R\rightarrow \infty,$$
where we have used 
the estimate $\|u\|_{L^2_t L^\frac{2n}{n-2}_x}<\infty$
(whose proof can be found in
\cite{bpst}).
\\
\\
{\em Proof of \eqref{due}} 
\\
\\

By using the H\"older inequality we get:
\begin{equation}\label{largevar}\left | \int \int_{|x|>1} 
\frac{|u|^2} {|x|^{3}} W\left( \frac x{|x|}\right) 
|\partial_{|x|} \phi_R |\hbox{ } dxdt\right |
\end{equation}
$$\leq  \|W\|_{L^\infty_x}\left [ \int ( 
\int |u(t)|^\frac{2n}{n-2} dx )^\frac{n-2}{n} dt \right]
\left (\int_{|x|>1} |x|^{-\frac {3n}2}|\partial_{|x|} 
\phi_R|^\frac n2  dx \right )^\frac 2n.$$

Notice also that $$\int_{|x|>1} 
|x|^{-\frac {3n}2}|\partial_{|x|} \phi_R|^\frac n2  dx$$
$$=\int_{1<|x|<R} |x|^{-\frac {3n}2}|\partial_{|x|} \phi_R|^\frac n2  dx
+ \int_{|x|>R} |x|^{-\frac {3n}2}|\partial_{|x|} \phi_R|^\frac n2  dx$$
$$\leq \frac 1{R^\frac n2} \int_{1<|x|<R} \frac 1{|x|^n} \hbox{ } dx
+ C  \int_{|x|>R} |x|^{-\frac {3n}2} dx$$
where we have used
$$\partial_{|x|} \phi\left(\frac xR \right)= 0\left( \frac{|x|}R\right).$$

By combining this estimate with
\eqref{largevar} and recalling that 
$\|u\|_{L^2_t  L^\frac{2n}{n-2}_x}<\infty$,
we get
$$\left | \int \int_{|x|>1} 
\frac{|u|^2} {|x|^{3}} W\left( \frac x{|x|}\right) 
|\partial_{|x|} \phi_R |\hbox{ } dxdt\right |
\rightarrow 0 \hbox{ as }
R\rightarrow \infty.$$

Next we treat the integral in the cylinder $\{|x|<1\}\times \mathbf R$.
We use the Taylor formula as above and we 
get:
$$\partial_{|x|}\phi(|x|)= 0(|x|) \hbox{ as } |x|\rightarrow 0$$ 
and then
$$\left | \int \int_{|x|<1} 
\frac{|u|^2} {|x|^{3}} W\left( \frac x{|x|}\right) 
|\partial_{|x|} \phi_R|\hbox{ } dxdt\right |
$$
$$
\leq  \frac CR
\left | \int \int_{|x|<1} 
\frac{|u|^2} {|x|^2} W\left( \frac x{|x|}\right) 
dxdt\right |\rightarrow 0 \hbox{ as } 
R\rightarrow \infty,
$$
where at the last step we have used \eqref{ulTI}. 

\hfill$\Box$

We are now able to prove
theorems \ref{mainSL}
and \ref{mainUCSL}.

\vspace{0.3cm}

\noindent {\bf Proof of theorems \ref{mainSL} and \ref{mainUCSL}}.

Due to a density argument it is sufficient 
to prove the theorems for $f \in H^1_x\cap L^2_{|x|^2}$.

Next we split the proof in two steps.

\vspace{0.1cm}

\noindent {\em Proof of r.h.s. in \eqref{eqSL}}

\vspace{0.2cm}

It is sufficient to replace in \eqref{BarcRuiVeg}
the generic function $\psi$
with the family of rescaled functions $R \phi\left (\frac xR \right),$
where $\phi$ is a function that satisfies proposition
\ref{brvpsi} and by estimating the r.h.s. in \eqref{BarcRuiVeg}
by using \eqref{dele}.\\

Let us point - out that the l.h.s.
in \eqref{eqSL}
follows from theorem \ref{mainUCSL}.
\\
\\
{\em Proof of theorem \ref{mainUCSL}} 
\\
\\

First of all let us notice that if we choose in the identity \eqref{BarcRuiVeg}
the function $\psi$ to be equal to the function $\phi$ given 
in proposition \ref{brvpsi}, then it is not difficult to verify
that  
\begin{equation*}
\int_{|x|>1} \frac{|\nabla_\tau u|^2}{|x|} \hbox{ } dx <\infty,
\end{equation*}
provided that $f\in H^1_x$,
and in particular
\begin{equation}\label{revised}
\lim_{R\rightarrow \infty}\int_{|x|>R} \frac{|\nabla_\tau u|^2}{|x|} \hbox{ } dx=0.
\end{equation}

Let us fix a function
$h(r)\in C^\infty_0(\mathbf R; [0, 1])$ such that:
$$h(r)\equiv 1 \hbox{  } \forall r\in {\mathbf R}
\hbox{ s.t. } |r|<1, h(r)\equiv 0 \hbox{  } \forall r \in {\mathbf R}
\hbox{ s.t. } |r|> 2, 
$$
$$h(r)= h(-r)\hbox{  }  \forall r \in \mathbf R.$$

We introduce also the functions
$\phi, H \in C^\infty( {\mathbf R})$:
$$\phi(r) = \int_0^r (r-s) h(s) ds \hbox{  } \hbox{ and } \hbox{  }
H(r)=\int_0^r h(s) ds,$$
(since now on in the proof the function $\phi$ will be the one 
defined above and not the one given in proposition \ref{brvpsi}).
Notice that 
\begin{equation}\label{differ}
\phi''(r)=h(r), \phi'(r)= H(r) \hbox{  } 
\forall r\in {\mathbf R} \hbox{ and }
\lim_{r\rightarrow \infty} \partial_{|x|} 
\phi(r)=\int_0^\infty h(s)ds.
\end{equation}

Moreover an elementary computation 
shows that: 

\begin{equation}\label{bilaplac}\Delta^2 \phi = 
\frac C{|x|^3} \hbox{  } \forall x\in 
{\mathbf R}^n \hbox{ s.t. } |x| \geq 2,
\end{equation}
where $\Delta^2$ is the bilaplacian operator.

Thus the function
$\phi$ defined above satisfies the assumptions of lemma \ref{noma}.
Notice also that the assumptions in proposition \ref{radiationSL}
are satisfied by $\phi$.

\vspace{0.1cm}

In the sequel we shall need the rescaled functions
$$\phi_{R}=R \phi\left (\frac{x}{R} \right ) \hbox{  }
\forall x\in {\mathbf R}^n
\hbox{ and }R>0,$$
(where $\phi$ is the function defined above) 
and we shall exploit the following 
elementary identity:
\begin{equation}\label{radiahess}\nabla \bar u D^2 \psi \nabla u=
\partial_{|x|}^2\psi |\partial_{|x|} u|^2 + 
\frac{\partial_{|x|} \psi}{|x|}|\nabla_\tau u|^2,\end{equation}
where $\psi$ is any regular radial function
and $u$ is another regular function.

By combining this identity with \eqref{BarcRuiVeg} 
and with proposition \ref{radiation2}, 
where we choose $\psi= 
\phi_{R}$, and recalling \eqref{differ} we get:
\begin{equation}\label{limitSLNL}
\int_0^\infty \int
\left [ \partial_{|x|}^2 \phi_{R} |\partial_{|x|} u|^2 
+ \frac{\partial_{|x|} \phi_{R}}
{|x|} | \nabla_\tau u|^2\right .
\end{equation}
\begin{equation*}
\left . -\frac 14
|u|^2 \Delta^2 \phi_{R}
+\frac{|u|^2} {|x|^{3}} W\left( \frac x{|x|}\right) \partial_{|x|} \phi_R 
\right ]dx dt 
\end{equation*}
\begin{equation*}
\geq  \frac 14 \left (\int_0^\infty h(s) ds
\right ) \|f \|_{\dot H^\frac 12_x }^2 
+\frac 12 {\mathcal Im} \int \bar f \hbox{ }
\nabla f \cdot \nabla \phi
\left( \frac x R \right ) dx
\hbox{  } \forall R>0.
\end{equation*}

By using now \eqref{uno} and \eqref{due} 
we get
\begin{equation*}
\lim_{R\rightarrow \infty} \int_0^\infty \int
\left [
- \frac 14 |u(t)|^2 \Delta^2 \phi_{R}  
+\frac{|u(t)|^2} {|x|^{3}} W\left( \frac x{|x|}\right) 
|\partial_{|x|} \phi_R |
\right ]dxdt=0.
\end{equation*}

On the other hand since $\phi$ is a radially symmetric function
we have $$\lim_{x\rightarrow 0} \nabla \phi(x)=0,$$
that due to the dominated convergence theorem implies
$$\lim_{R\rightarrow 0}
\int \bar f \hbox{ } \nabla f \cdot \nabla \phi
\left( \frac x R \right ) dx=0.$$

By combining these facts 
with \eqref{limitSLNL} we get
\begin{equation}\label{crusiaSL}
\liminf_{R\rightarrow \infty} 
\int_0^\infty\int
\left ( \partial_{|x|}^2\phi_{R} |\partial_{|x|} u|^2 
+ \frac{\partial_{|x|} \phi_{R}}
{|x|} | \nabla_\tau u|^2\right )
dx dt \end{equation}
\begin{equation*}\geq \frac 14 \left (\int_0^\infty h(s) ds
\right ) \|f \|_{\dot H^\frac 12_x}^2.
\end{equation*}

On the other hand by using 
the cut - off properties of $h$, \eqref{revised}
and noticing that $\partial_{|x|}\phi =0(|x|)$ as $|x|\rightarrow 0$, we get:
\begin{equation}\label{limsup} 
\liminf_{R\rightarrow \infty} \frac 1R \int_0^\infty \int_{B_{2R}}
(|\partial_{|x|} u|^2 + |\nabla_\tau u|^2) \hbox{ } dxdt 
\end{equation}
$$\geq \liminf_{R\rightarrow \infty} \int_0^\infty \int 
\left (\partial_{|x|}^2\phi_{R}|\partial_{|x|} u|^2 
+ \frac{\partial_{|x|} \phi_{R}}
{|x|} | \nabla_\tau u|^2\right ) dxdt 
$$$$\geq \frac 14 
\left (\int_0^\infty h(s) ds \right ) \|f\|_{\dot H^\frac 12_x}^2$$
where we have used \eqref{crusiaSL} at the last step.
The proof is complete.

\hfill$\Box$
 
\section{On the asymptotic behaviour
of solutions to critical NLS\\
and proof of theorems \ref{uniqueSN}, \ref{radiation}}

The main aim of this section
is to prove theorem \ref{radiation}
that represents the nonlinear version
of proposition \ref{radiationSL}.

First of all let us recall a precise statement about
the global existence result
to \eqref{SN} with small initial data.

\begin{thm}\label{basicNLS}
There exists $\epsilon_0>0$ such that for any
$f\in L^2_x$ 
with $\|f\|_{L^2_x}<\epsilon_0$,
the Cauchy problem
\eqref{SN} has an unique global solution 
$$u\in {\mathcal C}_t (L^2_x)\cap L^{2+\frac 4n}_{t,x}.$$
Moreover there exists a function
$$(0, \epsilon_0) \ni \epsilon 
\rightarrow R(\epsilon)\in {\mathbf R}^+$$
such that: 
\begin{enumerate}
\item 
$\lim_{\epsilon \rightarrow 0} R(\epsilon)=0$
\item  if $\|f\|_{L^2_x}<\epsilon$
then then unique solution to \eqref{SN} satisfies 
$\|u\|_{L^{2+\frac 4n}_{t,x}} < R(\epsilon)$.
\end{enumerate}
If moreover we assume $f\in H^1_x$, then
$u\in {\mathcal C}_t(H^1_x)$.
\end{thm}

\begin{remark}\label{SMA}
Notice that theorem \ref{basicNLS} provides a
global existence result
for small initial data, and also the arbitrary smallness of the 
$L^{2+\frac 4n}_{t,x}$ - norm 
of the solutions, provided  that the initial data are small enough.   
Let us recall that in the defocusing case it is sufficient
to assume $f\in H^1_x$ to have a global solution and no extra smallness
assumption is needed.
\end{remark}

In the sequel we shall make 
extensively use of the following inequality:
\begin{equation}\label{runstsick}
\| |D|^s (u|u|^q)\|_{L^p_x}
\leq C \||D|^s u\|_{L^r_x} \|u\|_{L^{sq}_x}^q
\end{equation}
where $C=C(s, q, p, r, s)>0$, $0\leq s \leq 1$,
$1<p, r, s<\infty$ and 
$$\frac 1p= \frac 1r + \frac 1s.$$

We shall also make use of the 
classical Strichartz inequalities
that we mention below for completeness (for a proof see \cite{kt}).

Assume that $u$ solves the following Schr\"odinger 
equation with forcing term:
$${\bf i}\partial_t u \pm \Delta u=F(t, x)$$
$$u(0)=f, (t, x)\in {\mathbf R}_t \times {\mathbf R}^n_x, n\geq 3$$
then the following a - priori estimates are satisfied:
$$\|u\|_{L^p_tL^q_x}\leq C\left (\|f\|_{L^2_x}+ 
\|F\|_{L^{\tilde p'}_tL^{\tilde q'}_x}\right ),$$
where
$$\frac{2}{p}+ \frac nq=\frac{2}{\tilde p}+ 
\frac n{\tilde q}=\frac n2, \hbox{  } 
2\leq p, \tilde p \leq \infty,$$
$$\frac 1{\tilde p}+ \frac 1{\tilde p'}=
\frac 1{\tilde q}+ \frac 1{\tilde q'}=1,$$
and $C=C(p, q, \tilde p, \tilde q)>0$.

\vspace{0.1cm}

Next we shall prove some preliminary propositions that 
will be useful along the proof
of theorem \ref{radiation}.

\begin{prop}\label{regNLS}
Let $v\in 
{\mathcal C}_t(H^1_x) \cap L_{t,x}^{2+\frac 4n}$ be the unique solution to 
$${\bf i}\partial_t v+ \Delta v \mp v|v|^\frac 4n=0, $$$$v(0)=g$$
where
$g\in L^2_{|x|}\cap H^1_x$
and $\|g\|_{L^2_x}<\epsilon$, with $\epsilon>0$ small enough.
Then the following estimate holds: 
\begin{equation}\label{impoNLS}
\|e^{-{\bf i}\frac{|x|^2}{4}} v(1)\|_{\dot H^\frac 12_x}
\leq  C \|g\|_{L^2_{|x|}}.
\end{equation}
\end{prop}

\noindent {\bf Proof.}
Let us introduce the
function
$$v^*(t, x)= \frac 1{t^\frac n2}e^{-{\bf i}\frac{|x|^2}{4t}}
v\left (\frac xt, \frac 1t \right).
$$

It is easy to verify
that $v^*$ satisfies:
\begin{equation}\label{vstar}{\bf i} \partial_t v^* - \Delta v^*
\pm v^* | v^*|^\frac 4n=0,
\end{equation} \begin{equation*}
v^*(1) = e^{-{\bf i}\frac {|x|^2}4} v(1)
\end{equation*}
and then
\begin{equation}\label{vstar11}{\bf i} \partial_t 
(|D|^s v^*) - \Delta (|D|^s v^*)
\pm |D|^s (v^* | v^*|^{\frac 4n})=0,
\end{equation}
\begin{equation*}
(|D|^s v^*)(1)= |D|^s(e^{-{\bf i}\frac {|x|^2}4} v(1)).
\end{equation*}

Hence $|D|^s v^*$ satisfies the free Schr\"odinger equation with
forcing term given by $\mp |D|^s (v^* | v^*|^{\frac 4n}$).

We can then apply Strichartz estimates 
for a suitable choice of the parameters $p, q, \tilde p, \tilde q$
in order to get:
$$\||D|^s v^* \|_{L^2_tL^{\frac{2n}{n-2}}_x}
\leq C \left (\||D|^s v^*(1)\|_{L^2_x} + \||D|^s (v^* | v^*|^{\frac 4n})
\|_{L^\frac{2n+4}{n+6}_tL^\frac{2n(n+2)}{n^2-4+4n}_x}  \right)
$$
$$
\leq C \left (\||D|^s v^*(1)\|_{L^2_x} + 
\||D|^s v^*\|_{L^2_tL^\frac{2n}{n-2}_x} 
\| v^*
\|_{L^{2+\frac 4n}_{t,x}}^{\frac 4n} \right)
$$
where at the last step we have used \eqref{runstsick}.
On the other hand 
$\|v^*\|_{L^{2+\frac 4n}_{t,x}}$ is small 
(this fact follows from the smallness of $v$ that in turn follows
from the smallness assumption done on $g$, 
see remark \ref{SMA}) and then we get
from the previous estimate the following one:
\begin{equation}\label{striend}
\||D|^s v^*\|_{L^2_t L^\frac{2n}{n-2}_x}
\leq C \||D|^s v^*(1)\|_{L^2_x}.
\end{equation}

Let us introduce now 
the functions $w(t, x)$ and $z(t, x)$ defined as the unique solutions to:
\begin{equation}\label{hom}
{\bf i}\partial_t z + \Delta z=0, 
\end{equation}
$$
z(0)=g$$
and
\begin{equation}\label{inhom}
{\bf i}\partial_t w+ \Delta w\mp v|v|^\frac 4n=0, 
\end{equation}$$
w(0)=0.
$$

It is clear that the following identity holds:
\begin{equation}\label{ideult}
z+w=v.
\end{equation}

Along with $z$ and $w$
we introduce also
$$w^* = 
\frac 1{t^\frac n2}e^{- {\bf i}\frac{|x|^2}{4t}} 
w\left (\frac 1t, \frac xt \right)$$
and
$$z^* = 
\frac 1{t^\frac n2}e^{-{\bf i}\frac{|x|^2}{4t}} 
z\left (\frac 1t, \frac xt \right).$$

Notice that due to \eqref{ideult} we have also
\begin{equation}\label{idebanal}
v^*=z^* + w^*.
\end{equation}

Next we shall estimate separately $z$ and $w$.

\vspace{0.1cm}

\noindent {\em Estimate for $z$}

\vspace{0.2cm}

Since $z$ is defined by \eqref{hom},
we can apply corollary \ref{pseudoNLS} for $t=1$, in order to deduce
\begin{equation}\label{partial1}
\||D|^\frac 12 z^*(1)\|_{L^2_x}=\||D|^\frac 12 [e^{-{\bf i} 
\frac{|x|^2}{4}} z(1)]\|_{L^2_x}\leq 
 \frac 1 {\sqrt 2} \|g\|_{L^2_{|x|}}.
\end{equation}

\vspace{0.1cm}

\noindent {\em Estimate for $w$}

\vspace{0.2cm}

Let us notice that the following identity trivially holds:
$$\|w^*(t)\|_{L^2_x}=\left \|w\left (\frac 1t\right )\right \|_{L^2_x}.$$

On the other hand (by definition) $w(0)=0$ and then 
due to the previous identity we get 
\begin{equation}\label{L2prime}\lim_{t\rightarrow \infty} 
\|w^* (t)\|_{L^2_x}=0.
\end{equation}

Moreover the function
$w^*$ satisfies: 
\begin{equation}\label{s}
{\bf i}\partial_t w^* - \Delta 
w^* = \mp v^* |v^*|^\frac 4n, 
\end{equation}
\begin{equation*}
w^*(1)= e^{-{\bf i}\frac{|x|^2}4}w(1)
\end{equation*}
and then 
\begin{equation}\label{s12}
{\bf i}\partial_t (|D|w^*) - \Delta 
(|D| w^*) = \mp |D|(v^* |v^*|^\frac 4n), 
\end{equation}
\begin{equation*}
|D|w^*(1)= |D|(e^{-{\bf i}\frac{|x|^2}4}w(1)).
\end{equation*}

We can combine again as above the
Strichartz estimates with \eqref{runstsick} in order to deduce:
\begin{equation}\label{second111}
\||D| w^* \|_{L^\infty_t L^2_x}
\leq C \left (\||D| w^*(1)\|_{L^2_x}+ 
\| |D|(v^*|v^*|^\frac 4n) \|_{L^\frac{2n+4}{n+6}_t
L^\frac{2n(n+2)}{n^2-4+4n}_x}
\right ) 
\end{equation}
\begin{equation*}
\leq C  \left( \||D| w^*(1)\|_{L^2_x} +
\||D| v^*\|_{L^2_t L^\frac{2n}{n-2}_x} 
\|v^*\|_{L^{2+\frac 4n}_{t,x}}^\frac 4n\right )
\end{equation*}
\begin{equation*}
\leq C \left (\||D| w^*(1)\|_{L^2_x}
+ \||D| v^*(1)\|_{L^2_x}\|v^*\|_{L^{2+\frac 4n}_{t,x}}^\frac 4n\right ),
\end{equation*}
where we have used \eqref{striend} for $s=1$
at the last step.

In particular we have
$$\sup_{t\in \mathbf R} \||D| w^*(t)\|_{L^2_x}<\infty.$$

By combining this fact with
\eqref{L2prime} and with the following elementary estimate:
$$\||D|^\frac 12 \phi \|_{L^2_x}
\leq C \|\phi\|_{L^2_x}^\frac 12
\| |D|\phi\|^\frac 12_{L^2_x} \hbox{   } \forall \phi\in H^1_x,$$
we deduce that
$$\lim_{t\rightarrow \infty} \||D|^\frac 12 w^*(t)\|_{L^2_x}=0.$$

As a consequence of this fact and \eqref{s}
we have that
\begin{equation}\label{s123}
{\bf i}\partial_t (|D|^\frac 12 w^*) - \Delta 
(|D|^\frac 12  w^*) = \mp |D|^\frac 12(v^* |v^*|^\frac 4n), 
\end{equation}
\begin{equation*}
(|D|^\frac 12 w^*)(\infty)=0.
\end{equation*}

By combining again Strichartz estimates  
(with the initial condition at infinity) with \eqref{runstsick}
we get: 
\begin{equation}\label{second}
\||D|^\frac 12 w^* \|_{L^\infty_t L^2_x}
\leq C
\| |D|^\frac 12 (v^*|v^*|^\frac 4n) \|_
{L^\frac{2n+4}{n+6}_tL^\frac{2n(n+2)}{n^2-4+4n}_x} 
\end{equation}
\begin{equation*}
\leq C  
\||D|^\frac 12 v^*\|_{L^2_tL^\frac{2n}{n-2}_x} 
\|v^*\|_{L^{2+\frac 4n}_{t,x}}^\frac 4n
\leq C \||D|^\frac 12 v^*(1)\|_{L^2_x}\|v^*\|_{L^{2+\frac 4n}_{t,x}}^\frac 4n,
\end{equation*}
where we have used \eqref{striend}
for $s=\frac 12$ at the last step.
\\
\\

By combining 
\eqref{idebanal}, \eqref{partial1} and \eqref{second} we get:
\begin{equation*}
\||D|^\frac 12 v^*(1)\|_{L^2_x}=\||D|^\frac 12 (z^*(1) + w^*(1))\|_{L^2_x}
\end{equation*}
\begin{equation*}
\leq 
C \left (
\||D|^\frac 12 v^*(1)\|_{L^2_x}\|v^*\|_{L^{2+\frac 4n}_{t,x}}^\frac 4n
+ \|g\|_{L^2_{|x|}}\right ),
\end{equation*}
that
due to the smallness of 
$\|v^*\|_{L^{2+\frac 4n}_{t,x}}$ implies:
\begin{equation*}
\|e^{-{\bf i}\frac{|x|^2}4} v(1)\|_{\dot H^\frac 12_x}
=\||D|^\frac 12 v^* (1)\|_{L^2_x} \leq C \|g\|_{L^2_{|x|}}. 
\end{equation*}

\hfill$\Box$

\begin{prop}
Let $u \in {\mathcal C}_t(H^1_x) \cap L_{t,x}^{2+\frac 4n}$ 
be the unique solution to 
$${\bf i}\partial_t u - \Delta u \pm u|u|^\frac 4n=0, $$$$u(0)=f$$
where $f\in H^1_x$
and $\|f\|_{L^2_x}<\epsilon$, with $\epsilon>0$ small enough.
Then the following estimate holds: 
\begin{equation}\label{giuste}
\|f\|_{\dot H^\frac 12_x} \leq C \|u(1)\|_{\dot H^\frac 12_x}
\end{equation}
where $C>0$ is a constant that does not depend on $f$.
\end{prop}

\noindent {\bf Proof.}
Notice that $w=|D|^\frac 12 u$ satisfies:
$${\bf i}\partial_t 
w
-\Delta w \pm |D|^\frac 12 (u|u|^\frac 4n)=0$$
$$w(1)=|D|^\frac 12 u(1),$$
then we can combine Strichartz inequalities with 
\eqref{runstsick} in order to get:
$$\|w\|_{L^2_tL^\frac{2n}{n-2}_x}
\leq C\left ( \||D|^\frac 12 u(1)\|_{L^2_x}+
\||D|^\frac 12 (u|u|^\frac 4n) 
\|_{L^\frac{2n+4}{n+6}_tL^\frac{2n(n+2)}{n^2-4+4n}_x}\right )
$$
$$\leq
C\left ( \||D|^\frac 12 u(1)\|_{L^2_x}+
\||D|^\frac 12 u\|_{L^2_tL^\frac{2n}{n-2}_x} \||u|^\frac 4n
\|_{L^\frac{n+2}2_{t,x}}\right )$$ 
$$=C\left ( \||D|^\frac 12 u(1)\|_{L^2_x}+
\||D|^\frac 12 u\|_{L^2_t L^\frac{2n}{n-2}_x} 
\|u\|^\frac 4n_{L^{2+\frac 4n}_{t,x}}\right ).$$ 
 
Due to the smallness of $\|u\|_{L^{2+\frac 4n}_{t,x}}$
(that comes from the smallness assumption done on $f$) we get
\begin{equation}\label{giustefra}
\||D|^\frac 12 u\|_{L^2_t L^\frac{2n}{n-2}_x}=
\|w\|_{L^2_t L^\frac{2n}{n-2}_x}
\leq 
C \||D|^\frac 12 u(1)\|_{L^2_x}.
\end{equation}

By using again Strichartz estimates 
with a different choice of the 
parameters $p, q$ we get
$$\|w\|_{L^\infty_t L^2_x}
\leq C\left (\||D|^\frac 12 u(1)\|_{L^2_x}
+ \| |D|^\frac 12 (u|u|^\frac 4n)
\|_{L^\frac{2n+4}{n+6}_tL^\frac{2n(n+2)}{n^2-4+4n}_x}\right)$$
$$\leq C\left (\||D|^\frac 12 u(1)\|_{L^2_x}
+ \||D|^\frac 12 u\|_{L^2_t L^\frac{2n}{n-2}_x}
\||u|^\frac 4n\|_{L^\frac{n+2}2_{t,x}} \right )$$
that is equivalent to
$$ \||D|^\frac 12 u\|_{L^\infty_t L^2_x}
\leq 
 C\left (\||D|^\frac 12 u(1)\|_{L^2_x}
+ \||D|^\frac 12 u\|_{L^2_t L^\frac{2n}{n-2}_x}
\|u\|_{L^{2+\frac 4n}_{t,x}}^\frac 4n \right ).$$

By combining this estimate with \eqref{giustefra} we get
$$\||D|^\frac 12 u\|_{L^\infty_t L^2_x}
\leq 
 C\||D|^\frac 12 u(1)\|_{L^2_x}
$$
and in particular
$$\|f\|_{\dot H^\frac 12_x}
= \|u(0)\|_{\dot H^\frac 12_x}
\leq 
\||D|^\frac 12 u\|_{L^\infty_t L^2_x}
$$$$\leq 
C \||D|^\frac 12 u(1)\|_{L^2_x}
= \|u(1)\|_{\dot H^\frac 12_x}.
$$

The proof is complete.

\hfill$\Box$

\vspace{0.2cm}

\noindent {\bf Proof of theorem \ref{radiation}.}
Let us recall that \eqref{conservation2} 
implies the following identity:
\begin{equation}\label{pseudoestimateSNL}
\left \|\frac xt u(t) - 2{\bf i}\nabla u(t)\right \|_{L^2_{x}}^2 
\pm \frac{2n}{n+2} \int |u(t)|^{2+\frac 4n } dx 
= \frac 1{t^2}\|f\|_{L^2_{|x|^2}}^2. 
\end{equation}

Since $u\in {\mathcal C}_t(H^1_x)$ it is meaningful to consider the trace
$u(\bar t)\in H^1_x$ for any $\bar t\in \mathbf R$ and since 
by assumption 
$u\in L^{2+\frac 4n}_{t,x}$,
we can deduce that
there exists a sequence $t_k\rightarrow \infty$
such that:
\begin{equation}
\int |u(t_k)|^{2+\frac 4n} dx\rightarrow 0 \hbox{ as } k\rightarrow \infty.
\end{equation}

By combining this fact with \eqref{pseudoestimateSNL} we get:
\begin{equation}\label{somSNL}
\lim_{k\rightarrow \infty}\left \|\frac x {t_k} u(t_k) - 2{\bf i} 
\nabla u (t_k)\right \|_{L^2_x}^2
=0.
\end{equation}

We split now the proof in two parts.
\\
\\
{\em Construction of $g$ and proof of \eqref{dotNLS}}
\\
\\

Let us recall that the conformal 
transfomation of $u(t, x)$:
$$\tilde u (t, x)= \frac{1}{t^\frac n2} e^{\frac{{\bf i}|x|^2}{4t}} 
u \left (\frac 1t, \frac xt \right),
$$
satisfies the Cauchy problem \eqref{SNP}
under the initial condition
$$\tilde u (1)= e^{\frac{{\bf i}|x|^2}{4}} u \left (1\right)$$
and in particular 
\begin{equation}\label{small}\|\tilde u (1)\|_{L^2_x}=\|u(1)\|_{L^2_x}
= \|f\|_{L^2_x}<\epsilon,
\end{equation}
where we have used the conservation of 
the charge for the unique solution to
\eqref{SN} (see \eqref{charge}).

Then $\tilde u$ satisfies the following Cauchy problem 
$$
{\bf i}\partial_t \tilde u + \Delta \tilde u\mp \tilde u 
|\tilde u|^\frac 4n=0,
$$$$\tilde u(1)\in L^2_x 
\hbox{ and } \|\tilde u(1)\|_{L^2_x}<\epsilon,
(t, x) \in (0, \infty) \times {\mathbf R}^n.$$

Due to the global well - posedness
of this Cauchy problem (see theorem \ref{basicNLS})  
we deduce that $\tilde u$ can be extended
as a solution 
to the whole space ${\mathbf R} \times {\mathbf R}^n$.

Moreover this extension belongs
to the functional space ${\mathcal C}_t(L^2_x)\cap L^{2+\frac 4n}_{t,x}$ and
in particular it is
well defined one unique
$g\in L^2_x$ such that
the following limit exists:
\begin{equation}\label{limNL}
\lim_{t\rightarrow 0} \tilde u (t, x)= g \in L^2_x.
\end{equation}
  
Hence $\tilde u(t, x)$ satisfies the following Cauchy problem
\begin{equation}
\label{SNg}{\bf i}\partial_t \tilde u + \Delta \tilde u  
\mp \tilde u |\tilde u|^\frac 4n=0,
\end{equation}
\begin{equation*}
\tilde u(0)=g.
\end{equation*}

Due to
\eqref{limNL} we can deduce that
$$
\lim_{t\rightarrow 0}\left \| u\left ( \frac 1t, x\right)- 
t^\frac n2 e^{-{\bf i} t \frac{|x|^2}{4}} 
g \left(tx\right)\right \|_{L^2_x}=\lim_{t\rightarrow 0}
\left \|\frac 1{t^\frac n2} e^{{\bf i} \frac{|x|^2}{4t}} 
u\left (\frac 1t, \frac xt \right)
- g \left (x\right )\right \|_{L^2_x}
=0,$$
and in particular
\begin{equation}\label{L2}
\lim_{t\rightarrow \infty} \left \|u(t) - 
\frac 1{t^\frac n2}e^{-{\bf i} \frac{|x|^2}{4t}} g\left
  (\frac{x}{t}\right )\right \|_{L^2_x}=0.
\end{equation}

On the other hand, since $\tilde u(t, x)$ satisfies \eqref{SNg},
we can apply \eqref{impoNLS}
in order to get:
$$\|e^{-{\bf i}\frac{|x|^2}4} \tilde u(1)\|_{\dot H^\frac 12_x}^2
\leq C \int |x||g(x)|^2 dx,$$
that due to the definition of $\tilde u (t, x)$ is equivalent to
$$\|u(1)\|_{\dot H^\frac 12_x}^2
\leq C \int |x||g(x)|^2 dx.$$

We can then combine this estimate with \eqref{giuste}
in order to get \eqref{dotNLS}.
\\  
\\
{\em Proof of \eqref{somcon7NLS}}
\\
\\

Due to \eqref{somSNL} we can deduce that
\begin{equation}\label{som2}\lim_{k\rightarrow \infty}
\left (\int \bar u(t_k) 
\nabla u(t_k) \cdot  \nabla \psi \hbox{ } dx 
+\frac {\bf i} {2t_k} \int |x| | u(t_k) |^2 \partial_{|x|}\psi 
\hbox{ } dx\right )=0
\end{equation}
and then
\begin{equation}\label{trunc}\limsup_{t\rightarrow \infty}
\left( -{\mathcal Im}
\int \bar u(t) \nabla u(t) \cdot 
\nabla \psi \hbox{ } dx \right)
\end{equation}
\begin{equation*}
\geq \frac 12
\liminf_{k\rightarrow \infty} \int |x| | u(t_k) |^2 
\partial_{|x|}\psi \hbox{ } \frac {dx}{t_k}.
\end{equation*}

Let us fix now a real number $R>0$ and let us notice that due to
the non - negativity assumption done on $\partial_{|x|}\psi$ we get:
\begin{equation}\label{limint}\int |x| | u (t_k)|^2 
\partial_{|x|}\psi \hbox{ } \frac{dx}{t_k}
\geq  \int_{|x|\leq Rt_k} |x| | u(t_k) |^2 
\partial_{|x|}\psi \hbox{ } \frac{dx}{t_k}
\end{equation}
$$=  \int_{|x|\leq Rt_k} |x| \left (|u(t_k)|^2 -\frac 1{t_k^n} 
\left | g\left (\frac x{t_k} \right )\right |^2\right )
\partial_{|x|}\psi\hbox{ }
\frac{dx}{t_k}
$$$$+ \int_{|x|\leq Rt_k} |x| \left [(\partial_{|x|} \psi 
-\psi'(\infty)\right]  
\left | g\left ( \frac x {t_k}\right)\right |^2 \frac{dx}{t_k^{n+1}}
$$
$$
+ \psi'(\infty) \int_{|x|\leq Rt_k} 
|x| \left |g\left (\frac x{t_k} \right )\right |^2 \frac{dx}{t_k^{n+1}}
$$
where $g$ is the function constructed in the previous step.

Notice that the following estimate is trivial:
\begin{equation}\label{lim4}
\int_{|x|\leq Rt_k}
|x| \left (|u(t_k)|^2 -\frac 1{t_k^n} 
\left | g\left (\frac x{t_k} \right )\right |^2\right )
\partial_{|x|}\psi \hbox{ } \frac{dx}{t_k}
\end{equation}
\begin{equation*}
\leq R \|\partial_{|x|}\psi\|_{L^\infty_x}\int_{|x|\leq Rt_k}
\left ||u(t_k)|^2 -\frac 1{t_k^n} 
\left | g\left (\frac x{t_k} \right )\right |^2\right |dx
\rightarrow 0 \hbox{ as } k\rightarrow \infty,
\end{equation*}
where in the last step we have used \eqref{L2}.

On the other hand due to the change of variable formula we can prove that:
$$\left |\int_{|x|\leq Rt_k} |x| \left [(\partial_{|x|} \psi 
-\psi'(\infty)\right] \left | g\left ( \frac x {t_k}\right)
\right |^2 \frac{dx}{t_k^{n+1}}\right |
$$
$$\leq R  \int_{|x|\leq R} \left |\partial_{|x|} \psi ( t_k x )
-\psi'(\infty)\right| |g(x)|^2 dx,$$ 
that in conjunction  with 
the dominated convergence theorem 
and with assumption \eqref{infbil} implies:
\begin{equation}
\label{lim7}\lim_{k\rightarrow \infty}\int_{|x|\leq Rt_k} 
|x| \left [(\partial_{|x|} \psi 
-\psi'(\infty)\right] 
\left | g\left ( \frac x {t_k}\right)\right |^2 \frac{dx}{t_k^{n+1}}
=0.\end{equation}

Due again to the change of variable formula we get
$$
\int_{|x| \leq Rt_k} |x| \left |g\left (\frac x{t_k} \right )\right |^2 
\frac{dx}{t_k^{n+1}}
= \int_{|x|\leq R } |x| |g(x)|^2 dx,
$$
and in particular
\begin{equation}\label{lim14}
\lim_{k\rightarrow \infty}
\psi'(\infty) \int_{|x|\leq Rt_k} 
|x| \left |g\left (\frac x{t_k} \right )\right |^2 \frac{dx}{t_k^{n+1}}
= \psi'(\infty)\int_{|x|\leq R} |x| |g(x)|^2 dx.
\end{equation}

By combining 
\eqref{lim4},\eqref{lim7}, \eqref{lim14} and
\eqref{limint} we can deduce:
\begin{equation}\label{limfin}
\liminf_{k\rightarrow \infty}\int 
|x| |u(t_k)|^2 \partial_{|x|}\psi \hbox{  } \frac{dx}{t_k}
\end{equation}
$$\geq \psi'(\infty) \int_{|x|\leq R} |x| |g(x)|^2 dx \hbox{  }
\forall R>0.
$$

Since $R>0$ is arbitrary, we can combine 
\eqref{trunc} with \eqref{limfin},
in order to deduce 
\eqref{somcon7NLS}.

\hfill$\Box$


\vspace{0.1cm}

\noindent {\bf Proof of theorem \ref{uniqueSN}}
The proof is similar to the one of theorem \ref{uniqueSL}.
Let $g$ denotes the function constructed in
theorem \ref{radiation}, then we shall show that:
\begin{equation}\label{uniqueness}
\liminf_{t\rightarrow \infty}\int 
|x| |u(t)|^2\frac{dx}{t}
\geq \int_{|x| \leq  R} |x| |g(x)|^2 dx \hbox{  }
\forall R>0.
\end{equation}

In fact we have:
\begin{equation}\label{limintSL3} \int 
|x|| u(t) |^2 \frac {dx}t
\geq \int_{|x|\leq Rt} |x| | u(t) |^2 \frac{dx}t
\end{equation}
$$=  \int_{|x|\leq Rt} |x| \left [|u(t)|^2 -\frac 1{t^n} 
\left | g\left (\frac xt \right )\right |^2\right ]
\frac{dx}t
+ \int_{|x|\leq Rt} 
|x| \left |g\left (\frac xt \right )\right |^2 \frac{dx}{t^{n+1}} \hbox{ }
\forall R>0.
$$

Notice that the following estimate is trivial:
\begin{equation}\label{lim4SL3}
\int_{|x|\leq Rt}
|x| \left (|u(t)|^2 -\frac 1{t^n} 
\left | g\left (\frac xt \right )\right |^2\right )\frac{dx}t
\end{equation}
\begin{equation*}
\leq R \int_{|x|\leq Rt}
\left ||u(t)|^2 -\frac 1{t^n} 
\left | g\left (\frac xt \right )\right |^2\right |dx
\rightarrow 0 \hbox{ as } t\rightarrow \infty,
\end{equation*}
where in the last step we have used \eqref{L2}.

Due again to the change of variable formula we get
\begin{equation}\label{lasseb3}
\int_{|x| \leq Rt} |x| \left |g\left (\frac xt \right )\right |^2 \frac{dx}{t^{n+1}}
= \int_{|x|\leq R } |x| |g(x)|^2 dx.
\end{equation}

Hence \eqref{uniqueness} follows by combining 
\eqref{limintSL3}, \eqref{lim4SL3},
\eqref{lasseb3}.

In particular if 
$$\liminf_{t\rightarrow \infty}\int 
|x| |u(t)|^2 \hbox{ } \frac{dx}t=0,$$ then  
\eqref{uniqueness} implies $g\equiv 0$, and in turn  due to 
\eqref{L2} it gives $\lim_{t\rightarrow \infty} \|u(t)\|_{L^2_x}=0$.
By combining this fact with the conservation 
of the charge \eqref{charge}
we get $f\equiv 0$ and hence $u\equiv 0$.

\hfill$\Box$

\section{Proof of theorem \ref{mainSN}}

We shall follow basically  the strategy that has been used
in the proof
of theorems \ref{mainSL} and \ref{mainUCSL}.
More precisely we multiply
\eqref{SN} by the quantity  given in \eqref{multiplier},
we integrate on the strip $(0, T)\times {\mathbf R}^n$
and with elementary computations we get: 
\begin{equation}\label{BarcRuiVegNLS}
\int_{0}^T\int \left[\nabla \bar u  D^2 \psi  \nabla u
-\frac 14  |u|^2\Delta^2 \psi 
\pm \frac{1}{n+2}|u|^{2+\frac 4n} \Delta \psi  
\right ]dx dt\end{equation}
\begin{equation*} = -\frac 12
{\mathcal Im}
\int \bar u(T) \nabla u(T) \cdot \nabla \psi\hbox{ } dx
+\frac 12 {\mathcal Im}\int \bar f \hbox{ }
\nabla f \cdot \nabla \psi\hbox{ } dx.
\end{equation*}

The following lemma will be very important in order to prove the r.h.s. 
estimate in \eqref{eqSN}.
\begin{lem}
Assume that $u\in {\mathcal C}_t(H^1_x)\cap L^{2+\frac 4n}_{t,x}$
solves \eqref{SN} with $f\in H^1_x$ and
$\|f\|_{L^2_x}<\epsilon$, where $\epsilon>0$ is a 
small number. Assume moreover that
$\phi$ satisfies the same assumption as in proposition
\ref{brvpsi},
then there exists a constant $C>0$ that does not depend on $R>0$
and such that following estimates hold:
\begin{equation}\label{radNLSle}
\left |\int u(t) \nabla \bar u(t) \cdot \nabla \phi_R \hbox{ } dx
\right |\leq C \|f\|_{\dot H^\frac 12_x}^2 
\hbox{ } \forall t \in \mathbf R, R>0
\end{equation}
and
\begin{equation}\label{snls}
\int \int 
 |u|^{2+\frac 4n} |\Delta \phi_R |
\hbox{ } dxdt \leq C 
R^{\frac 2n -1 }\epsilon^2 \|f\|_{\dot H^\frac 12_x}^{\frac 4n} 
\hbox{  } \forall R>0,
\end{equation}
where
$\phi_R= R \phi\left(\frac xR \right)$.
\end{lem}

\noindent {\bf Proof.}

\vspace{0.1cm}

\noindent {\em Proof of \eqref{radNLSle}}

\vspace{0.2cm}

The proof of \eqref{radNLSle} is similar to the proof of \eqref{dele} 
provided that we are able to show the following a - priori bound:
\begin{equation}\label{auxbil}\|u(t)\|_{\dot H^\frac 12_x} \leq C 
\|f\|_{\dot H^\frac 12_x}
\hbox{ } \forall t \in \mathbf R,\end{equation}
where $u$ and $f$ are as in the statement.

In order to prove this inequality we notice that:
$${\bf i}\partial_t (|D|^\frac 12 u )
- \Delta (|D|^\frac 12 u)=\mp |D|^\frac 12(u|u|^\frac 4n ), 
$$$$
|D|^\frac 12 u(0)=|D|^\frac 12 f $$
hence by combining Strichartz estimates
with \eqref{runstsick}, and following the proof of
\eqref{striend}, 
we can get the following estimate:
\begin{equation}\label{striendp}
\||D|^\frac 12u \|_{L^2_t L^\frac{2n}{n-2}_x}
\leq C \|f\|_{\dot H^\frac 12_x}.
\end{equation}

On the other hand by combining again Strichartz estimates
with \eqref{runstsick} we get:
\begin{equation}\label{secondend}
\||D|^\frac 12  u\|_{L^\infty_t L^2_x}
\leq C \left (\||D|^\frac 12 f\|_{L^2_x} +
\| |D|^\frac 12 (u|u|^\frac 4n) 
\|_{L^\frac{2n+4}{n+6}_tL^\frac{2n(n+2)}{n^2-4+4n}_x}\right ) 
\end{equation}
\begin{equation*}
\leq C \left (\||D|^\frac 12 f\|_{L^2_x}+ 
\||D|^\frac 12 u\|_{L^2_tL^\frac{2n}{n-2}_x} 
\|u\|_{L^{2+\frac 4n}_{t,x}}^\frac 4n\right )
\end{equation*}
\begin{equation*}
\leq C \left (\||D|^\frac 12  f\|_{L^2_x}+\||D|^\frac 12 f
\|_{L^2_x}\|u\|_{L^{2+\frac 4n}_{t,x}}^\frac 4n\right ),
\end{equation*}
where we have used \eqref{striendp} at the last step.

Hence we can deduce \eqref{auxbil} by using the boundedness
of $\|u\|_{L^{2+\frac 4n}_{t,x}}$,
that in turn comes from the smallness assumption done on $f$.

\vspace{0.1cm}

\noindent {\em Proof of \eqref{snls}}

\vspace{0.2cm}

Since
$\phi$ satisfies proposition \ref{brvpsi},
it is easy to deduce that 
\begin{equation}\label{pointlap}|\Delta \phi|\leq \frac{C}{1+|x|} 
\hbox{  } \forall x\in {\mathbf R}^n.\end{equation}

Next we shall need the following a - priori estimates:
\begin{equation}\label{stri2le}
\|u\|_{L^2_tL^\frac{2n}{n-2}_x}\leq C \|f\|_{L^2_x}
\end{equation}
and 
\begin{equation}\label{stri1le}
\|u\|_{L^\infty_t \dot H^\frac 12_x}\leq C \|f\|_{\dot H^\frac 12_x}.
\end{equation}

Notice that \eqref{stri1le} is equivalent 
to \eqref{auxbil}, then we shall show \eqref{stri2le}.

Since $u$ solves \eqref{SN} we are in position to apply
the Strichartz estimates
in order to deduce:
$$\|u\|_{L^2_tL^\frac{2n}{n-2}_x}
\leq C\left ( \|f\|_{L^2_x}+
\|u|u|^\frac 4n\|_{L^\frac{2n+4}{n+6}_tL^\frac{2n(n+2)}{n^2-4+4n}_x}\right )
$$
$$\leq
C\left ( \|f\|_{L^2_x}+
\|u\|_{L^2_tL^\frac{2n}{n-2}_x} \||u|^\frac 4n
\|_{L^\frac{n+2}2_{t,x}}\right )$$ 
$$=C\left ( \|f\|_{L^2_x}+
\|u\|_{L^2_tL^\frac{2n}{n-2}_x} \|u
\|^\frac 4n_{L^{2+\frac 4n}_{t,x}}\right ),$$ 
that due to the smallness of 
$ \|u\|^\frac 4n_{L^{2+\frac 4n}_{t,x}}$
(that depends as usual on the smallness assumption
done on $f$) implies \eqref{stri2le}.

The H\"older inequality
implies the following chain of estimates: 
\begin{equation}\label{prert2}\int \int 
|u|^{2+\frac 4n} |\Delta \phi_R |
\hbox{ } dxdt = \frac 1 R \int \int 
|u|^{2+\frac 4n} \left |\Delta \phi \left (\frac xR\right )\right |
dxdt\end{equation}
$$\leq \frac 1R \left (
\int \left |\Delta \phi\left (\frac xR \right)
\right |^\frac{n^2}2 dx\right )^\frac 2{n^2}
\int \|u(t)\|_{L^{\frac{2n(n+2)}{n^2-2}}_x}^{2+\frac 4n} dt
$$
$$\leq 
\frac 1R \left (\int \left |\Delta 
\phi\left (\frac xR \right)\right |^\frac{n^2}2 dx \right )^\frac 2{n^2}
\int \|u(t)\|_{L^\frac{2n}{n-1}_x}^{\theta (2+\frac 4n)} 
\|u(t)\|_{L^\frac{2n}{n-2}_x}^{(1-\theta) (2+\frac 4n)}dt,$$
where
$$\frac{n^2-2}{2n(n+2)}=\frac {\theta(n-1)} {2n}
+ \frac{(1-\theta)(n-2)}{2n},$$
or equivalently $\theta= \frac 2{n+2}$. 

Then the previous estimate becomes
\begin{equation*}\int \int 
|u|^{2+\frac 4n} |\Delta \phi_R |
\hbox{ } dxdt \end{equation*}
\begin{equation*}
\leq \frac 1R \left (\int 
\left |\Delta \phi\left (\frac xR \right)
\right |^\frac{n^2}2 dx \right )^\frac 2{n^2}
\int \|u(t)\|_{L^\frac{2n}{n-1}_x}^{\frac 4n} 
\|u(t)\|_{L^\frac{2n}{n-2}_x}^2dt
\end{equation*}
$$\leq \frac 1R \left (\int \left 
|\Delta \phi\left (\frac xR \right)
\right |^\frac{n^2}2 dx \right )^\frac 2{n^2}
\|u(t)\|_{L^\infty_t L^\frac{2n}{n-1}_x}^{\frac 4n} 
\|u(t)\|_{L^2_t L^\frac{2n}{n-2}_x}^2,$$
that due to the Sobolev embedding
\begin{equation}\label{sobolev}
\dot H^\frac 12_x \rightarrow L^\frac{2n}{n-1}_x
\end{equation}
implies:
\begin{equation}\label{inca}\int \int 
|u|^{2+\frac 4n} |\Delta \phi_R |
\hbox{ } dxdt \end{equation}
$$\leq \frac 1R \left (\int \left |\Delta \phi\left (\frac xR \right)
\right |^\frac{n^2}2 dx \right )^\frac 2{n^2}
\|u(t)\|_{L^\infty_t \dot H^\frac{1}{2}_x}^{\frac 4n} 
\|u(t)\|_{L^2_t L^\frac{2n}{n-2}_x}^2.$$

Notice that due to \eqref{pointlap} we get 
$$\left \|\Delta \phi \left( \frac xR\right)\right \|_{L_x^\frac{n^2}{2}}
\leq C R^\frac 2n 
\left (\int \frac{1}{(1+|x|)^{\frac{n^2}{2}}} \hbox{ } dx\right)^\frac{2}{n^2}
\hbox{ } \forall R>0.$$

By combining now this estimate with \eqref{stri2le},
\eqref{stri1le},
and \eqref{inca}
we deduce that
\begin{equation}\label{inca1}\int \int 
|u|^{2+\frac 4n} |\Delta \phi_R |
\hbox{ } dxdt \end{equation}
\begin{equation*}\leq CR^{\frac 2n-1} 
\|f\|_{\dot H^\frac 12_x}^\frac 4n \|f\|_{L^2_x}^2
\leq C R^{\frac 2n -1}\epsilon^2\|f\|_{\dot H^\frac 12_x}^\frac 4n 
\hbox{ } \forall R>0.
\end{equation*}

\hfill$\Box$

\vspace{0.2cm}

We shall need also the following

\begin{lem}\label{nomaNLS}
Let 
$u\in {\mathcal C}_t(L^2_x)\cap L^{2+\frac 4n}_{t,x}$ 
be the unique solution
to \eqref{SN} where $\|f\|_{L^2_x}<\epsilon$,
with $\epsilon>0$ small, then:
\begin{equation}\label{unoNLS}
\lim_{R\rightarrow \infty}  \int \int
|u|^2 |\Delta^2 \phi_R| 
\hbox{ } dxdt =0\end{equation} 
and
\begin{equation}\label{dueNLS}
\lim_{R\rightarrow \infty}  \int \int 
|u|^{2+\frac 4n} |\Delta \phi_R |\hbox{ } dxdt=0,
\end{equation}
where $\phi\in C^\infty({\mathbf R}^n)$ is a radially symmetric function
such that 
$$
|\partial_{|x|}\phi|\leq C, |\Delta^2 \phi |\leq \frac{C}
{(1+|x|)^3} \hbox{  } \forall x\in {\mathbf R}^n$$
and $\phi_R=R \phi\left (\frac{x}R\right )$.
\end{lem}

\noindent{\bf Proof.}
\\
\\
\noindent {\em Proof of \eqref{unoNLS}}
\\
\\

Let us notice that since $u$ solves \eqref{SN}
we can apply Strichartz estimates in order to get:
\begin{equation}\label{endpointNLS}\|u\|_{L^2_t L^\frac{2n}{n-2}_x}
\leq C\left( \|f\|_{L^2_x} + \||u|^{1+\frac 4n} 
\|_{L^\frac {2(n+2)}{n+4}_{t,x}}\right)
\end{equation}
$$=C \left ( \|f\|_{L^2_x} + \|u\|_{L^{2+\frac 4n}_{t,x}}^{1+\frac 4n}
\right)<\infty.
$$ 

On the other hand the H\"older inequality implies:
$$\int \int |u|^2 |\Delta^2 \phi_R(x)| \hbox{ } dx dt
\leq \left [\int \left ( \int |u|^\frac{2n}{n-2} 
\hbox{ } dx \right )^\frac{n-2}{n} dt\right ]
\left ( \int |\Delta^2 \phi_R(x)|^\frac n2  \hbox{ } dx \right )^\frac 2n$$
$$\leq  C \|u\|_{L^2_t L^\frac{2n}{n-2}_x}^2 
\left (\int  \frac 1{(R + |x|)^\frac {3n}2} \hbox{ } dx \right)^\frac 2n
\rightarrow 0 \hbox{ as
} R\rightarrow \infty,$$
where we have used in the last step
the fact that $\|u\|_{L^2_t  L^\frac{2n}{n-2}_x}<\infty$
that comes from \eqref{endpointNLS}.
\\
\\
{\em Proof of \eqref{dueNLS}} 
\\
\\

By using the H\"older inequality we get:
$$\int \int 
|u|^{2+\frac 4n} |\Delta \phi_R |\hbox{ } dxdt
\leq  \|\Delta \phi_R\|_{L^\infty_x} 
\int \int |u|^{2+\frac 4n} \hbox{ } dx dt$$
$$=\frac 1R \|\Delta \phi\|_{L^\infty_x} 
\|u\|_{L^{2+\frac 4n}_{t,x}}^{2+\frac 4n}
\rightarrow 0 \hbox{ as } R \rightarrow \infty.$$

\hfill$\Box$

\noindent {\bf Proof of theorem \ref{mainSN}}
\\
\\

Let us first prove the estimate
\begin{equation}\label{luva}\sup_{R>R_0}\frac 1R 
\int_{0}^\infty \int_{|x|<R}
|\nabla u|^2 \hbox{ } dxdt\end{equation}
$$\leq C \left( \|f\|_{\dot H^\frac 12_x}^2 + \frac 1{R_0^{1-\frac 2n}}
\|f\|_{\dot H^\frac 12_x}^\frac 4n
\right ) \hbox{ } \forall R_0>0.$$

Notice that \eqref{BarcRuiVegNLS} implies:
\begin{equation*}
\int_{0}^T\int \left[\nabla \bar u D^2 \psi \nabla u
-\frac 14  |u |^2\Delta^2 \psi \right ]dx dt =
\end{equation*}$$\mp \frac 1{n+2}\int^T_{0} \int
|u|^{2+\frac 4n} \Delta \psi \hbox{ } 
dx dt 
$$$$- \frac 12
{\mathcal Im}
\int \bar u(T) \nabla u(T) \cdot \nabla \psi \hbox{ }dx
+ \frac 12 {\mathcal Im}\int \bar f \hbox{ }
\nabla f\cdot \nabla \psi \hbox{ } dx.
$$

If we choose in the previous identity
$\psi=R \phi\left( \frac xR\right )$ (with $R>R_0$)
where $\phi$ is as in proposition \ref{brvpsi},
and if we recall \eqref{radNLSle},\eqref{snls}, 
then we can deduce \eqref{luva}.
Notice that \eqref{luva} trivially implies
$$\limsup_{R\rightarrow \infty}
\frac 1R \int \int_{B_R}|\nabla u|^2 \hbox{ } dxdt
\leq C \|f\|_{\dot H^\frac 12_x}^2.$$

Next we shall prove
\begin{equation}\label{zar}\liminf_{R\rightarrow \infty}
\frac 1R \int \int_{B_R}|\nabla u|^2 \hbox{ } dxdt
\geq c \|f\|_{\dot H^\frac 12_x}^2\end{equation}
and it will be sufficient to complete
the proof of the theorem.

In fact the proof of \eqref{zar} can be done by exploiting the identity
\eqref{BarcRuiVegNLS}, where we choose 
the function $\psi$ to be equal to the
functions $\phi_R$ used in the proof of theorem \ref{mainUCSL}.

Then the argument follows as in the proof 
of theorem \ref{mainUCSL} with some minor changes.
In fact it is sufficient to replace proposition \ref{radiationSL}
with theorem \ref{radiation},
and to use \eqref{unoNLS} and \eqref{dueNLS} 
instead of \eqref{uno} and \eqref{due}.

\hfill$\Box$

\section{Appendix}

In order to make this paper self - contained
we shall give in this appendix a result already presented in \cite{BRV}. 
More precisely we shall prove the existence of a test function $\phi$ 
with suitable properties
that has been extensively used along this paper. 
\begin{prop}\label{brvpsi}
Assume that $n\geq 3$, then there exists
a radially symmetric function $\phi:{\mathbf R}^n\rightarrow \mathbf R$
such that:
\begin{enumerate}
\item $\phi(0)=\partial_{|x|}\phi(0)=0$ and $\partial_{|x|}^2\phi(0)> 0$;
\item for any $x\in {\mathbf R}^n$ we have
\begin{equation*}
\Delta^2 \phi  \leq 0 \hbox{ } \forall x\in {\mathbf R}^n;
\end{equation*}
\item $\partial_{|x|} \phi, \partial_{|x|}^2 \phi > 0
\hbox{ } \forall x \in {\mathbf R}^n\setminus \{0\};$
\item there exists $C>0$ such that
$$\partial_{|x|}\phi, |x|\partial_{|x|}^2\phi \leq C
\hbox{  } \forall x \in {\mathbf R}^n;
$$
\item the following limit exists 
$$\lim_{|x|\rightarrow \infty}\partial_{|x|}\phi\in (0, \infty).$$
\end{enumerate}
In particular we have:
\begin{equation}\label{12tor}|\Delta \phi|\leq \frac{C}{1+|x|}\hbox{ } \forall x\in {\mathbf R}^n,
\end{equation}
\begin{equation}\label{13tor}
\nabla u D^2 \phi \nabla \bar u \geq C (|\partial_{|x|}^2 u|^2
+ |\nabla_{\tau} u|^2)\hbox{ } \forall x\in {\mathbf R}^n
\hbox{ s.t. } |x|<1,
\end{equation}
and 
\begin{equation}\label{15tor}
\nabla u D^2 \phi \nabla \bar u \geq 
C \frac{|\nabla_{\tau} u|^2}{|x|}\hbox{ } \forall x\in {\mathbf R}^n
\hbox{ s.t. } |x|>1,
\end{equation}
for any function $u$.
\end{prop}

\noindent {\bf Proof.}

It is easy to show that \eqref{12tor} follows by 
writing the laplacian in polar coordinates,
while \eqref{13tor} and \eqref{15tor} follow from the identity
\eqref{radiahess}.

Next we shall focus on the construction of $\phi$ that satisfies $(1), (2), (3), (4)$.

The main strategy is to solve the following equation:
$$-\Delta^2 \phi= h_\eta(|x|)$$
where \begin{equation}\label{heta}
h_\eta (|x|)=\chi_{\{|x|<1\}}+ 
\frac \eta{|x|^3}\chi_{\{|x|>1\}}
\end{equation}
and $\eta\geq 0$ is a suitable parameter that will be choosen later.
 
The equation $-\Delta^2 \phi=h_\eta$ can be written in polar coordinates
in the following equivalent way:
$$-r^{-(n-1)} \partial_r (r^{n-1} 
\partial_r (r^{-(n-1)} \partial_r (r^{n-1} 
\partial_r \phi(r))))=h_\eta(r).$$  

By integrating directly this equation we get:
\begin{equation}\label{ode} \partial_r \phi=-\frac 1{r^{n-1}}
\int_0^r u^{n-1}\int_0^u \frac 1{s^{n-1}}
\int_0^s t^{n-1} h(t)dt + \lambda r,
\end{equation}
where $\lambda\in \mathbf R$ is a generic 
number that will be choosen later.

Next we split the proof in two cases.

\vspace{0.1cm}

\noindent {\em First case: $n=3$}

\vspace{0.2cm}

Let us choose $h_\eta$ as in \eqref{heta} with $\eta=0$, 
in this way by an explicit 
integration we get:
$$\partial_r \phi = \lambda r - \frac{r^3}{30}
\hbox{ } \forall \hbox{ } 0<r<1$$
$$\partial_r \phi = \lambda r +\frac 16 - \frac r 6 
- \frac 1{30 r^2}\hbox{ } \forall \hbox{ } r>1.$$

In particular if we choose $\lambda= \frac 16$ then we deduce
that $\partial_r \phi$ satisfies the desired assumptions.
Moreover with the previous choice of $\lambda$ we have:
$$\partial_r^2 \phi = \frac 16 - \frac{r^2}{10}
\hbox{ } \forall \hbox{ } 0<r<1$$
$$\partial_r^2 \phi = \frac 1{15 r^3}\hbox{ } \forall \hbox{ } r>1.$$

It is now easy to check that all the properties 
required to $\phi$ are fulfilled.

\vspace{0.1cm}

\noindent{\em Second case: $n>4$}

\vspace{0.2cm}

In this case we choose $\eta>0$ in \eqref{heta}.
The precise value of $\eta$ will be choosen later.

An explicit integration of \eqref{ode} gives:
$$\partial_{r} \phi=\lambda r - \frac{r^3}{2n(n+2)} 
\hbox{ } \forall \hbox{ } 0<r<1$$
$$\partial_{r} \phi 
=\lambda r + \frac \eta{(n-1)(n-3)}
- \frac{(1+2\eta)n-(3+6\eta)}{2n(n-2)(n-3)}r - \frac{\eta n -n +3}
{2n(n-2)(n-3)}\frac{1}{r^{n-3}}
$$$$- \left[ \frac{(1-\eta)n^4+(3\eta -6)n^3 +(11+4\eta)n^2-
(12 \eta +6)n }{2n^2(n-1)(n-2)(n-3)(n+2)}
\right ] \frac{1}{r^{n-1}}\hbox{ } \forall \hbox{ } r>1.$$

Notice that if we choose 
\begin{equation}\label{lamb}
\lambda=\frac{(1+2\eta)n -(3+6\eta)}{2n(n-2)(n-3)}
\end{equation}
then we get
$$\partial_{r}^2 \phi=\lambda  - \frac{3 r^2}{2n(n+2)} 
\hbox{ } \forall \hbox{ } 0<r<1$$
$$\partial_{r}^2 \phi 
=\frac{\eta n -n +3}{2n(n-2)}\frac{1}{r^{n-2}}
$$$$+ \left[ \frac{(1-\eta)n^4+(3\eta -6)n^3 + (11+4\eta)n^2
-(12 \eta +6)n }{2n^2(n-2)(n-3)(n+2)}
\right ] \frac{1}{r^{n}}\hbox{ } \forall \hbox{  } r>1.$$

It is easy to verify that all the properties required
to $\phi$ will be fulfilled provided that
we can choose $\lambda$ and $\eta$ that satisfy \eqref{lamb}
and moreover
$$\lambda>\frac{3}{2n(n+2)}, \hbox{ } \eta n -n +3>0,$$$$
(1-\eta)n^4+(3\eta -6)n^3 + (11+4\eta)n^2-(12 \eta +6)n >0,$$
$$\frac \eta{(n-1)(n-3)}>\frac{\eta n - n +3}{2n(n-2)(n-3)}
$$
$$+ \frac{(1-\eta)n^4+ (3\eta -6)n^3 + (11+4\eta)n^2
-(12 \eta +6)n }{2n^2(n-1)(n-2)(n-3)(n+2)}.$$

Explicit computations show that
the previous conditions are equivalent to look 
for a suitbale $\eta>0$ such that:
$$Max \left \{\frac{n-3}{n},\frac{-2n^2+8n - 6}{n^3-2n^2-5 n+6}
\right \} < \eta
<\frac{n^3-6n^2+11n-6}{(n-2)(n+2)(n-3)}.$$

An elementary computation shows that
$$\frac{-2n^2+8n - 6}{n^3-2n^2-5 n+6}
<0$$
for $n>3$, it is then enough to verify that
$$\frac{n-3}{n} < \frac{n^3-6n^2+11n-6}{(n-2)(n+2)(n-3)}.$$

It is easy to verify that this condition is fulfilled for any $n\geq4$.


\hfill$\Box$

\end{document}